\documentclass{amsart}
\usepackage{graphicx}
\usepackage[all]{xy}
\usepackage{latexsym}
\usepackage{amssymb}
\usepackage{amsmath}
\usepackage{color}

\usepackage{tikz, tikz-cd}
 \tikzstyle{int}=[circle, draw,fill=black,outer sep=0,minimum size=3pt, inner 
 \tikzstyle{ext}=[circle, draw=black,outer sep=0,inner sep=1pt]

\usepackage{fullpage}
\xyoption{arc}

 \newcommand{\ac}{\scriptstyle \textrm{!`}}
\def\id{{\mbox{1 \hskip -7pt 1}}}
\newcommand{\sgn}{{\mathit s  \mathit g\mathit  n}}
 \newcommand{\lon}{\longrightarrow}
 \newcommand{\bu}{\bullet}
 
 \newcommand{\rar}{\rightarrow}

\newcommand{\p}{{\partial}}
\newcommand{\Id}{{\mathrm{Id}}}

\newcommand{\Der}{\mathrm{Der}}

 \newcommand{\Z}{{\mathbb Z}}
 \newcommand{\bS}{{\mathbb S}}

 \newcommand{\K}{{\mathbb K}}

 \newcommand{\ot}{\otimes}

\newcommand{\Def}{\mathsf{Def}}

 \newcommand{\Beq}{\begin{equation}}
 \newcommand{\Eeq}{\end{equation}}
 \newcommand{\Beqr}{\begin{eqnarray}}
 \newcommand{\Eeqr}{\end{eqnarray}}
 \newcommand{\Beqrn}{\begin{eqnarray*}}
 \newcommand{\Eeqrn}{\end{eqnarray*}}
 \newcommand{\Ba}{\begin{array}}
 \newcommand{\Ea}{\end{array}}
 \newcommand{\Bi}{\begin{itemize}}
 \newcommand{\Ei}{\end{itemize}}
 \newcommand{\Bc}{\begin{center}}
 \newcommand{\Ec}{\end{center}}

 \newcommand{\cA}{{\mathcal A}}
 \newcommand{\cB}{{\mathcal B}}
 
 \newcommand{\caD}{{\mathcal D}}
 \newcommand{\cE}{{\mathcal E}}
 \newcommand{\cF}{{\mathcal F}}

 \newcommand{\cI}{{\mathcal I}}

 \newcommand{\cM}{{\mathcal M}}
 
 \newcommand{\cP}{{\mathcal P}}

 \newcommand{\cV}{{\mathcal V}}

 \newcommand{\ga}{\gamma}
 
 \newcommand{\Ga}{\Gamma}

 \newcommand{\la}{\lambda}

 \newcommand{\Hom}{{\mathrm H\mathrm o\mathrm m}}
 
 \newcommand{\sip}{\smallskip}
 \newcommand{\bip}{\bigskip}
 \newcommand{\mip}{\vspace{2.5mm}}

\newcommand{\wh}{\widehat}

\newcommand{\RGC}{\mathsf{RGC}}
 \newcommand{\Ass}{\mathcal{A}ss}

\newcommand{\ORGC}{\mathsf{ORGC}}
\newcommand{\dRGC}{\mathsf{dRGC}}

\newcommand{\Lie}{\mathcal{L} \mathit{ie}}

\newcommand{\rgc}{\mathfrak{rgc}}
\newcommand{\orgc}{\mathfrak{orgc}}

 \newcommand{\bbu}{\mbox{\resizebox{2.8mm}{!}{$\bullet$}}}

\newcommand{\PCY}{\mathcal{P}re\mathcal{CY}^3}
\newcommand{\PCYf}{\mathcal{P}re\mathcal{CY}}
\newcommand{\DPois}{\mathcal{DP}\mathit{ois}}

\theoremstyle{plain}
\swapnumbers

\newtheorem{prop-def}[theorem]{Proposition-definition}

\newtheorem{f-theorem}{Formality Theorem}[section]
\newtheorem{main-theorem}{Main~Theorem}[section]
\newtheorem{section-theorem}{Theorem}[section]

\theoremstyle{definition}

\begin{document}

 \sloppy

 \newenvironment{proo}{\begin{trivlist} \item{\sc {Proof.}}}
  {\hfill $\square$ \end{trivlist}}

\long\def\symbolfootnote[#1]#2{\begingroup%
\def\thefootnote{\fnsymbol{footnote}}\footnote[#1]{#2}\endgroup}

\title{On  dg properads of pre-CY algebras}

\author{Sergei A. Merkulov}
\address{Sergei~A.\ Merkulov,  Faculty of Mathematics, National Research University Higher School of Economics, Moscow }
\email{smerkulov25@gmail.com}

\begin{abstract} We present new small models for the dg properads that govern pre-CY algebras: one for general pre-CY algebras and one for pre-CY algebras without curvature terms. These small models have only four and, respectively, three generators with valencies $\leq 4$ (in contrast to the original properads which have infinitely many generators with arbitrarily large valencies). We prove that the valency upper bound  $4$ is sharp.

\sip

We study the cohomology group of the deformation complex of the dg properad
governing general  pre-CY algebras and prove that it contains the totality
$
\prod_{g\geq 1}H^\bu(\cM_{g,1})
$
of cohomology groups of moduli spaces $\cM_{g,1}$ of genus $g$ algebraic curves with one marked point. As an application of this result we prove that the Tradler-Zeinalian properad of $\cV^{(d)}$-algebras is not Koszul.

\sip

\noindent {\sc Mathematics Subject Classifications} (2000). 14H10, 18G85, 18M70

\noindent {\sc Key words}. Properad, pre-Calabi-Yau algebra, moduli spaces of algebraic curves.  
\end{abstract}

 \maketitle
\markboth{}{}

{\Large
\section{\bf Introduction}
}

Pre-CY algebras play a significant role in modern non-commutative geometry, topology,  and representation theory. They were introduced around 2013 by Maxim Kontsevich and Yannis Vlassopoulos (see \S 1.4 in \cite{KTV}); in  finite dimensions, they are equivalent to  $\cV^{(d)}_\infty$-algebras   studied earlier in \cite{TZ}. The fact that the pre-CY algebras are well-defined in infinite dimensions is of the  utmost importance for many applications.

\sip

For any (possibly, infinite-dimensional) graded vector space $A$ and any integer $d\in \Z$, there is an associated 
{\em higher}\,
 Hochschild graded vector space
$$
C_{[d]}^\bu(A) :=\prod_{m\geq 1} C_{[d]}^{(m)}(A):= \prod_{m\geq 1}
\left( \bigoplus_{n_1,\ldots, n_m\geq 0\atop m+n_1+...+n_m\geq 2}
\Hom\left(\bigotimes_{i=1}^{m}
(A[1])^{\ot n_i}, (A[2-d])^{\ot m}\right)\right)_{\Z_m},
$$
the space of (co)invariants under the natural diagonal action of the cyclic
 groups $\Z_m$ on the vector spaces
$\Hom\left(\bigotimes_{i=1}^{m} (A[1])^{\ot n_i}, (A[2-d])^{\ot m}\right)$. The space $C_{[d]}^\bu(A)$ has a natural
$\Lie_{d}$-algebra\footnote{If $\Lie$ stands for the standard operad of Lie algebras, then $\Lie_d:=\Lie\{d-1\}$ is by definition its degree shifted version which controls graded Lie algebras with Lie bracket of degree $1-d$.}  structure \cite{IKV}. Its Maurer-Cartan elements $\pi$, that is
 degree $d$ elements 
 $$
 \pi=\sum_{m\geq 1} \pi_m, \ \ \ \pi_m\in  C_{[d]}^{(m)}(A),
 $$
 satisfying the equation 
 $
[\pi,\pi]=0
$,
are called  {\em pre-CY algebra structures on $A$ of degree $d$} \cite{IKV, KTV}. The Lie algebra $C_{[d]}^\bu(A)$ containing the standard Hochschild complex as a Lie subalgebra  $C_{[d]}^{(1)}(A)$ so that the summand $\pi^1$
makes $A$ into an $\Ass_\infty$-algebra.

\sip

There is a dg properad $\PCYf_d=\{ \PCYf_d(m,n)\}_{m\geq 1, n\geq 0}$ whose representations in a graded vector space $A$ are in one-to-one correspondence with degree $d$ pre-CY algebra structures in $A$. This properad is best understood as spanned by {\it ribbon}\, graphs with labelled incoming and labelled outgoing hairs (or legs), e.g.
\Beq\label{1: exs of haired graphs}
\Ba{c}\resizebox{21mm}{!}{
\mbox{$\xy
 (0,16)*{^{\bar{1}}}="out1";
 (0,-3)*{_{\bar{2}}}="out2";
 (0,-18)*{_3}="in3";
 (16,4)*{^{1}}="in1";
  (16,-4)*{^{2}}="in2";
 (-8,0)*{\bbu}="0";
 (0,8)*{\bbu}="1";
  (8,0)*{\bbu}="2";
(0,-12)*{\bbu}="3";
 \ar @{<-} "1";"0" <0pt>
  \ar @{<-} "1";"2" <0pt>
   \ar @{->} "2";"3" <0pt>
 \ar @{->} "2";"0" <0pt>
  \ar @{<-} "0";"3" <0pt>
  \ar @{->} "3";"out2" <0pt>
  \ar @{->} "1";"out1" <0pt> 
   \ar @{->} "in3";"3" <0pt>
    \ar @{<-} "2";"in1" <0pt> 
     \ar @{<-} "2";"in2" <0pt> 
\endxy$}}
\Ea
\in  \PCYf_d(3,2),
\ \ \ \ 
\Ba{c}\resizebox{13mm}{!}{
\mbox{$\xy
 (0,16)*{^{\bar{1}}}="out1";
 (-8,0)*{\bbu}="0";
 (0,8)*{\bbu}="1";
  (8,0)*{\bbu}="2";
(0,-12)*{\bbu}="3";
 \ar @{<-} "1";"0" <0pt>
  \ar @{<-} "1";"2" <0pt>
   \ar @{<-} "2";"3" <0pt>
 \ar @{->} "2";"0" <0pt>
  \ar @{<-} "0";"3" <0pt>
  \ar @{->} "1";"out1" <0pt> 
\endxy$}}
\Ea
\in  \PCYf_d(1,0).
\Eeq
The properad $\PCYf_d$ is free,  generated by all ribbon\footnote{As usual, the adjective {
\it ribbon}\, means that all $m+n$ legs attached to the vertex $\bu$ of an $(m,n)$-corolla are cyclically ordered.} \ $(m,n)$-corollas with labelled $n\geq 0$ incoming legs  and $m\geq 1$ outgoing legs such that $m+n\geq 2$, e.g.
\Beq\label{1: exs of ribb corollas}
\underbrace{\Ba{c}\resizebox{11mm}{!}{\xy
 (-5,5)*{^{\bar{1}}}="a",
(0,0)*{\bullet}="b",
(5,5)*{^{\bar{2}}}="c",
\ar @{<-} "a";"b" <0pt>
\ar @{->} "b";"c" <0pt>
\endxy} \Ea}_{m=2,n=0}\ \ , \ \ \ \
\underbrace{\xy
(-5,-6)*{_{\bar{1}}}="a1",
(-5,6)*{^{\bar{2}}}="a2",
(5,5)*{^{\bar{3}}}="a3",
(5,-6)*{^{3}}="a4",
(8,0)*{\, ^2}="a5",
(-8,0)*{\, ^1}="a6",
(0,0)*{\bu}="b",
\ar @{<-} "a1";"b" <0pt>
\ar @{<-} "a2";"b" <0pt>
\ar @{<-} "a3";"b" <0pt>
\ar @{->} "a4";"b" <0pt>
\ar @{->} "a5";"b" <0pt>
\ar @{->} "a6";"b" <0pt>
\endxy}_{m=n=3}\ \ , 
\ \
\underbrace{\Ba{c}\resizebox{16mm}{!}{  \xy
(0,8)*{^{\bar{1}}}="U";
(-5,5)*{^{\bar{3}}}="UL";
(5,5)*{^{\bar{2}}}="UR";
(0,-7)*{_{\bar{4}}}="D";
    (0,0)*{\bu}="C";
(-7,1)*{^3}="L1";
(-5,-4)*{_2}="L2";    
  (6,-3)*{^1}="R";
  
\ar @{->} "C";"U" <0pt>
\ar @{->} "C";"UL" <0pt>
\ar @{->} "C";"UR" <0pt>
\ar @{->} "C";"D" <0pt>
\ar @{<-} "C";"L1" <0pt>
\ar @{<-} "C";"L2" <0pt>
\ar @{<-} "C";"R" <0pt>
 \endxy}
 \Ea}_{m=4,n=3}\ \in \PCYf_d.
\Eeq
Every ribbon $(m,n)$-corolla is assigned the cohomological degree $(2-d)m -n+d$. The differential $\delta$ in $\PCYf_d$ is given by splitting the vertices (see \S 2.1.1 for a precise definition). The dg properad $\PCYf_d$ can be identified with the cobar construction on the cooperad $(\cP\cV^{(1-d)})^{\ac}$ of the quadratic properad $\cP\cV^{(1-d)}$ of so called $\cV^{(d)}$-algebras introduced  in \cite{TZ}.

\sip

Our main purpose in this paper is to find the smallest possible quasi-isomorphic version of $\PCYf_d$. 

\sip

Let  $\langle J, \delta J\rangle$ be the differential closure of an ideal $J$ in $\PCYf_d$   generated by (i) all ribbon $(m,0)$-corollas with $m\geq 3$, and (ii) all  ribbon $(m,n)$-corollas with $m+n\geq 4$ except the following one 
\Beq\label{1: DLie corolla}
  \Ba{c}\resizebox{13mm}{!}{  \xy
(0,-5)*{}="U";
(0,+5)*{}="D";
    (0,0)*{\bu}="C";
  (5,0)*{}="R";
   (-5,0)*{}="L";
\ar @{->} "C";"U" <0pt>
\ar @{->} "C";"D" <0pt>
\ar @{<-} "C";"L" <0pt>
\ar @{<-} "C";"R" <0pt>
 \endxy}
 \Ea.
\Eeq
Let, and consider a quotient dg properad
$$
\Delta\PCYf_d:= \PCYf_d/\langle J, \delta J\rangle
$$
generated just by four ribbon corollas,
\Beq\label{1: types of corollas in DeltaPCYf}
\Ba{c}\resizebox{8mm}{!}{\xy
 (0,0)*{}="a",
(4,-3)*{\bullet}="b",
(8,0)*{}="c",
\ar @{<-} "a";"b" <0pt>
\ar @{->} "b";"c" <0pt>
\endxy} \Ea, 
 \ \ \
 \Ba{c}\resizebox{8mm}{!}{  \xy
(0,-6)*{}="1";
    (0,-0.2)*{\bu}="L";
  (-4,5)*{}="C";
   (4,5)*{}="D";
\ar @{<-} "D";"L" <0pt>
\ar @{<-} "C";"L" <0pt>
\ar @{->} "1";"L" <0pt>
 \endxy}
 \Ea, 
 \ \ \ 
 \Ba{c}\resizebox{8mm}{!}{  \xy
(0,6)*{}="1";
    (0,0.2)*{\bu}="L";
  (-4,-5)*{}="C";
   (+4,-5)*{}="D";
\ar @{->} "D";"L" <0pt>
\ar @{->} "C";"L" <0pt>
\ar @{<-} "1";"L" <0pt>
 \endxy}
 \Ea,
 \ \ \  
  \Ba{c}\resizebox{13mm}{!}{  \xy
(0,-5)*{}="U";
(0,+5)*{}="D";
    (0,0)*{\bu}="C";
  (5,0)*{}="R";
   (-5,0)*{}="L";
\ar @{->} "C";"U" <0pt>
\ar @{->} "C";"D" <0pt>
\ar @{<-} "C";"L" <0pt>
\ar @{<-} "C";"R" <0pt>
 \endxy}
 \Ea\ .
\Eeq
subject to quadratic relations described explicitly in \S 2.

\mip

\subsection{Main Theorem}\label{1: Theorem on PCYf} {\it  The epimorphism of dg properads 
$$
\pi: \PCYf_d\lon \Delta \PCYf_d
$$
is a quasi-isomorphism.}

\mip

It is worth noting that the above mentioned quadratic relations for the generators of $\Delta \PCYf_d$ are nice enough to guarantee that there is an 
epimorphism of dg properads
$$
p: \Delta \PCYf_d \lon \DPois_d \lon 0
$$
where $\DPois_d$ is the properad of (degree shifted) double Poisson algebras \cite{vdB} (equipped with the zero differential). The kernel of $p$ is the {\it differential}\, ideal in  $\Delta \PCYf_d$ generated by the following two ribbon corollas
$$
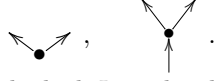
\Ba{c}\resizebox{8mm}{!}{\xy
 (0,0)*{}="a",
(4,-3)*{\bullet}="b",
(8,0)*{}="c",
\ar @{<-} "a";"b" <0pt>
\ar @{->} "b";"c" <0pt>
\endxy} \Ea, 
 \ \ \
 \Ba{c}\resizebox{8mm}{!}{  \xy
(0,-6)*{}="1";
    (0,-0.2)*{\bu}="L";
  (-4,5)*{}="C";
   (4,5)*{}="D";
\ar @{<-} "D";"L" <0pt>
\ar @{<-} "C";"L" <0pt>
\ar @{->} "1";"L" <0pt>
 \endxy}
 \Ea.
$$

One can consider next a differential ideal $I$ in the dg properad $\PCYf_d$ generated by all ribbon corollas of the form (``curvature terms")
  \Beq\label{1: sources}
\Ba{c}\resizebox{23mm}{!}{\xy
(2,-6)*{\ldots},
(-13,+7)*{},
(-8,-1)*{},
(-3,-7)*{},
(7,-7)*{},
(13,-7)*{},
 (0,0)*{\bu}="a",
(0,7)*{}="0",
(-7,+7)*+{_k}="b_1",
(-9,-3)*+{_1}="b_2",
(-3,-9)*+{_2}="b_3",
(8,-7)*+{}="b_4",
(5,9)*+{_{k-1}}="b_5",
%
\ar @{->} "a";"b_2" <0pt>
\ar @{->} "a";"b_3" <0pt>
\ar @{->} "a";"b_1" <0pt>
\ar @{->} "a";"b_4" <0pt>
\ar @{->} "a";"b_5" <0pt>
\endxy}\Ea, \ \ \ {k\geq 2},
\Eeq
and consider the quotient dg properad
$$
\PCY_d:= \PCYf_d/I,
$$
where the superscript $3$ emphasizes the fact that all generators of the free dg properad $\PCY_d$ are at least trivalent.
This properad was first studied by Alexandre Quesney in \cite{Q}, who proved that $\PCY_d$ can be identified with the cobar construction of the cooperad
$ \mathsf{BIB}_d^{\ac}$, which is Koszul dual to the so-called properad $ \mathsf{BIB}_d$ of {\it balanced infinitesimal bialgebras}.

\sip

Consider an ideal $J^3$ in the non-differential properad $\PCY_d$ generated by all possible ribbon $(m,n)$-corollas with $m+n\geq 4$ except the one 
in (\ref{1: DLie corolla}), and let $\langle J^3, \delta J^3 \rangle$ be its differential closure in $(\PCY_d, \delta)$.

\subsection{Theorem}\label{1: Theorem on PCY} {\it  The epimorphism of dg properads 
$$
\pi^3: \PCY_d\lon  \Delta\PCY_d:= \PCY_d/\langle J^3, \delta J^3 \rangle
$$
is a quasi-isomorphism.}

\mip

The dg properad  $\Delta\PCY_d$ has only three generators 
\Beq\label{1: generators of in DeltaPCY}
 \Ba{c}\resizebox{8mm}{!}{  \xy
(0,-6)*{}="1";
    (0,-0.2)*{\bu}="L";
  (-4,5)*{}="C";
   (4,5)*{}="D";
\ar @{<-} "D";"L" <0pt>
\ar @{<-} "C";"L" <0pt>
\ar @{->} "1";"L" <0pt>
 \endxy}
 \Ea, 
 \ \ \ 
 \Ba{c}\resizebox{8mm}{!}{  \xy
(0,6)*{}="1";
    (0,0.2)*{\bu}="L";
  (-4,-5)*{}="C";
   (+4,-5)*{}="D";
\ar @{->} "D";"L" <0pt>
\ar @{->} "C";"L" <0pt>
\ar @{<-} "1";"L" <0pt>
 \endxy}
 \Ea,
 \ \ \  
  \Ba{c}\resizebox{13mm}{!}{  \xy
(0,-5)*{}="U";
(0,+5)*{}="D";
    (0,0)*{\bu}="C";
  (5,0)*{}="R";
   (-5,0)*{}="L";
\ar @{->} "C";"U" <0pt>
\ar @{->} "C";"D" <0pt>
\ar @{<-} "C";"L" <0pt>
\ar @{<-} "C";"R" <0pt>
 \endxy}
 \Ea\ .
\Eeq
subject to quadratic relations described explicitly in \S 2 below.

\sip

We prove Theorems {\ref{1: Theorem on PCYf}} and {\ref{1: Theorem on PCY}}
in \S 2 below.

\subsection{On the sharpness of the conditions defining small models} A natural question: can one further enlarge the ideals $J$ and $J^3$ (and hence make smaller the associated quotient dg properads $\Delta\PCYf_d$ and $\Delta\PCY_d$) in such a way that Theorems {\ref{1: Theorem on PCYf}} and {\ref{1: Theorem on PCY}} hold true? Using the theory of cohomology groups of moduli spaces $\cM_{g,m}$ we show in \S 3 that the answer is {\bf no}.
This result follows from the following Theorem proved in \S 3 below.

\subsection{Theorem}\label{1: Theorem on Def complex} {\it Let $\Def(\PCYf_d \rar \PCYf_d)$ be the deformation complex of the identity autmorphism of the dg properad $\PCYf_d$. Then its cohomology $H^\bu(\Def(\PCYf_d \rar \PCYf_d))$ contains  contains the totality
$
\prod_{g\geq 1}H^\bu(\cM_{g,1})
$
of cohomology groups of moduli spaces $\cM_{g,1}$ of genus $g$ algebraic curves with one marked point.}

\sip

Let $V$ be the differential ideal in $(\PCYf_d,\delta)$  generated by all ribbon corollas except the following pair
\Beq\label{1: generators of PV(d)}
\Ba{c}\resizebox{8mm}{!}{\xy
 (0,0)*{}="a",
(4,-3)*{\bullet}="b",
(8,0)*{}="c",
\ar @{<-} "a";"b" <0pt>
\ar @{->} "b";"c" <0pt>
\endxy} \Ea, 
 \ \ \
  \Ba{c}\resizebox{8mm}{!}{  \xy
(0,6)*{}="1";
    (0,0.2)*{\bu}="L";
  (-4,-5)*{}="C";
   (+4,-5)*{}="D";
\ar @{->} "D";"L" <0pt>
\ar @{->} "C";"L" <0pt>
\ar @{<-} "1";"L" <0pt>
 \endxy}.
 \Ea
\Eeq 
 The quotient properad $\cP\cV^{(1-d)}$ defined by the short exact sequence
$$
0\lon  V \lon \PCYf_d\stackrel{j}\lon \cP\cV^{(1-d)}\lon 0
$$
is precisely the properad of $\cV^{(1-d)}$-algebras introduced in \cite{TZ}.  Theorem  {\ref{1: Theorem on Def complex}} implies that the epimorphism $j$ in the above short exact sequence can {\it not} a quasi-isomorphism which in turn implies the following

\subsubsection{\bf Corollary}\label{1: Corollary on TZ} {\it The Tradler-Zeinalian properad $\cP\cV^{(1-d)}$ is not Koszul}.

\sip

This corollary answers in negative Question 4.2 raised in \cite{PT}. 
This negative answer contrasts with the fact proven in \cite{PT} that the {\it dioperad}\,  $\cV^{(1-d)}$ underlying the properad  $\cP\cV^{(1-d)}$ {\it is}\, Koszul. A second proof of the Koszulness of the dioperad $\cV^{(1-d)}$ is given in \cite{Kh}.

\sip

Non-Koszulness of  Quesney's properad $ \mathsf{BIB}_d$ was established earlier in \cite{Me1} using a similar argument.

  \subsection{Some notation} We work over a field $\K$ of characteristic zero.
 The set $\{1,2, \ldots, n\}$ is denoted by $[n]$;  its group of automorphisms is
denoted by $\bS_n$; the trivial (resp., the sign) one-dimensional representation of
 $\bS_n$ is denoted by $\id_n$ (resp.,  $\sgn_n$). The cardinality of a finite set $S$ is denoted by $\# S$, while its linear span over a
field $\K$ by $\K\left\langle S\right\rangle$.


\sip

If $V=\oplus_{i\in \Z} V^i$ is a graded vector space, then
$V[k]$ stands for the graded vector space with $V[k]^i:=V^{i+k}$. For $v\in V^i$ we set $|v|:=i$.

\sip

{\bf Acknowledgements}. This work was supported by the research project HSE-BR-2025-84 implemented as a part of the Basic Research Program at the HSE University.

\bip

{\Large
\section{\bf Small models for dg properads $\PCYf_d$ and $\PCY_d$}
}

\sip

\subsection{Definition of $\PCYf_d$} A graph with a cyclic order of the half-edges incident to each vertex is called a {\it ribbon graph}.
A ribbon graph can be equivalently understood as a {\it fat graph}, that is,  a graph embedded in an orientable surface such that each vertex $v$ is a small disk $D$, and its edges are thin strips (ribbons)  glued along disjoint arcs $C_{edge}$ of the disk boundaries $\p D$. For each vertex-disk $v=D$ the complement $\p D\setminus \coprod C_{edge}$ is also a disjoint union of arcs 
called the set $C(v)$ of {\it corners}\, of the vertex $v$.
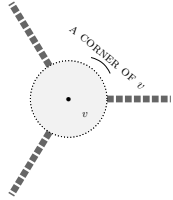
$$
\Ba{c}\resizebox{23mm}{!}{ 
\begin{tikzpicture}

    \draw[draw=black, thick, dotted, fill=black!5] (0,0) circle (1.2);


    \draw[line width=6pt, draw=black!60, dotted, rounded corners=3pt] (1.2,0) -- (3.2,0);
    \filldraw[fill=black!60, draw=black!60] (3.2,0) ;

    \draw[line width=6pt, draw=black!60, dotted, rounded corners=3pt] (-0.6,1.039) -- (-1.8,3.0);
    \filldraw[fill=black!60, draw=black!60] (-1.8,3.0);

    \draw[line width=6pt, draw=black!60, dotted, rounded corners=3pt] (-0.6,-1.039) -- (-1.8,-3.0);
    \filldraw[fill=black!60, draw=black!60] (-1.8,-3.0) ;

    \draw[draw=black, thick, dotted] (0,0) circle (1.2);

    \node at (0,0) {$\bullet$};
    \node[below right] at (0.3,-0.3) {$v$};
    
\draw[-, thick, bend right=25] (1.3,0.8) to node[above, sloped, font=\Large] {\sc a corner of $v$} (0.7,1.3);

\end{tikzpicture}
  }\Ea
   $$ 
   This set $C(v)$ is cyclically ordered. One can glue new thin strips to such corners, producing new ribbon graphs (an operation that  is used many times in this paper). A boundary $b$ of a ribbon graph consists of boundaries of thin strips and of corners of vertices. The set of corners belonging to a boundary $b$ is denoted by $C(b)$; this set is also cyclically ordered.

\sip

 Fix an arbitrary integer $d\in \Z$, and let $E_d(m,n)$ be the vector space spanned over a field $\K$ by all possible ribbon corollas with labelled $n\geq 0$ incoming legs  and $m\geq 1$ outgoing legs such that $m+n\geq 2$. By definition, each generator of the space $E_d(m,n)$ is assigned the cohomological degree $(2-d)m-n+d$. Here are some examples of ribbon corollas, 
\Beq\label{2: exs of ribb corollas}
\Ba{c}\resizebox{11mm}{!}{\xy
 (-5,5)*{^{\bar{1}}}="a",
(0,0)*{\bullet}="b",
(5,5)*{^{\bar{2}}}="c",
\ar @{<-} "a";"b" <0pt>
\ar @{->} "b";"c" <0pt>
\endxy} \Ea  \in E_d(2,0)
\ \ , \ \ \ \
\xy
(-5,-6)*{_{\bar{1}}}="a1",
(-5,6)*{^{\bar{2}}}="a2",
(5,5)*{^{\bar{3}}}="a3",
(5,-6)*{^{3}}="a4",
(8,0)*{\, ^2}="a5",
(-8,0)*{\, ^1}="a6",
(0,0)*{\bu}="b",
\ar @{<-} "a1";"b" <0pt>
\ar @{<-} "a2";"b" <0pt>
\ar @{<-} "a3";"b" <0pt>
\ar @{->} "a4";"b" <0pt>
\ar @{->} "a5";"b" <0pt>
\ar @{->} "a6";"b" <0pt>
\endxy\in E_d(3,3) 
\ \ , \ \
\Ba{c}\resizebox{16mm}{!}{  \xy
(0,8)*{^{\bar{1}}}="U";
(-5,5)*{^{\bar{3}}}="UL";
(5,5)*{^{\bar{2}}}="UR";
(0,-7)*{_{\bar{4}}}="D";
    (0,0)*{\bu}="C";
(-7,1)*{^3}="L1";
(-5,-4)*{_2}="L2";    
  (6,-3)*{^1}="R";
  
\ar @{->} "C";"U" <0pt>
\ar @{->} "C";"UL" <0pt>
\ar @{->} "C";"UR" <0pt>
\ar @{->} "C";"D" <0pt>
\ar @{<-} "C";"L1" <0pt>
\ar @{<-} "C";"L2" <0pt>
\ar @{<-} "C";"R" <0pt>
 \endxy}
 \Ea \in E_d(4,3).
\Eeq
The permutation group $\bS_m^{op}\times \bS_n$ acts on the vector space
$E_d(m,n)$ by relabelling outgoung and incoming hairs so that the collection
$$
E_d:=\{ E_d(m,n)\} 
$$
is an $\bS$-bimodule. Consider the free properad 
$$
\PCYf_d:= \cF ree \langle E_d\rangle
$$
generated by the $\bS$-bimodule $E_d$; properadic compositions in $\PCYf_d$ are given simply by gluing out-hairs to in-hairs as usual.  Hence the generators  of  $\PCYf_d$ can be understood as ribbon graphs $\Ga$ with labelled hairs and {\it directed}\, internal edges such that $\Ga$ has no closed paths of directed edges, see (\ref{1: exs of haired graphs}) for some examples. Let us denote by $E_{int}(\Ga)$  (resp. $V(\Ga)$) the set of internal edges (resp.\ of vertices) of such a hairy ribbon graph; the set of outgoing (resp. incoming) hairs of $\Ga$ is denoted by $H_{out}(\Ga)$
(resp., $H_{in}(\Ga)$).

\sip

Every hairy  ribbon graph $\Ga$ spanning $\PCYf_d$ comes equipped (often tacitly) with an orientation $or(\Ga)$,  defined as an ordering of {\it internal}\, edges and incoming hairs 
\Beq\label{2: or for d even}
or(\Ga):= \left(\bigwedge_{e \in E_{int}(\Ga)}e \right)\wedge \left(\bigwedge_{h \in H_{in}(\Ga)} h\right)  \ \ \  \text{for $d$ even},
\Eeq
 or an ordering of vertices and all hairs
 \Beq\label{2: or for d odd}
or(\Ga):= \left(\bigwedge_{v \in V(\Ga)} v\right)\wedge \left(\bigwedge_{\bar{h} \in H_{out}(\Ga)}\bar{h} \right)\wedge \left(\bigwedge_{h \in H_{in}(\Ga)} h\right) \ \ \  \text{for $d$ odd}.
\Eeq
There are two possible orientations, $or(\Ga)$ and $-or(\Ga)$ for each hairy ribbon graph $\Ga$ and we identify $(\Ga, or(\Ga))=-(\Ga,-or(\Ga)$, and abbreviate every such a pair $(\Ga, or(\Ga))$ to simply $\Ga$. Of course, the ordering of hairs can be chosen in agreement with their numerical labels, but this is not compulsory.   

\sip

The graded vector space of all possible representations of the non-differential properad $\PCYf_d$ in a graded vector space $A$
$$
\rho: \PCYf_d \lon \cE nd_A
$$
can be identified with the higher Hochschild space $C_{[d]}^\bu(A)$ of $A$ (see \S 1 for its definition)  as any such a representation $\rho$ is uniquely determined by its values on the ribbon  $(m,n)$-corollas which can be chosen arbitrarily.

\subsubsection{\bf The differential} The differential $\delta$ in the properad $\PCYf_d$ is defined via its action on the vertices of hairy ribbon graphs,
\Beq\label{2: delta in PreCYd}
\delta\Ga:= \sum_{v\in V(\Ga)} 
\delta_v\Ga,
\Eeq
and the graph $\delta_v\Ga$ is obtained from $\Ga$  by substituting into the vertex $v$ the ribbon graph $\xy
 (0,0)*{\bu}="a",
(5,0)*{\bu}="b",
\ar @{->} "a";"b" <0pt>
\endxy$ and then taking a sum over all possible reattachments of the edges (attached earlier to $v$) among the two newly created vertices in a way which respects cyclic order of edges and makes each new  vertex at least bivalent and having at least one outgoing hairs. 
Here are examples of such an action on 2, 3 and 4-valent vertices,
\Beq\label{2 d on 2 and 3val}
\delta_v 
\left(\Ba{c}\resizebox{9mm}{!}{ \xy
 (0,0)*{}="a",
(4,-3)*{\bullet}="b",
(8,0)*{}="c",
\ar @{<-} "a";"b" <0pt>
\ar @{->} "b";"c" <0pt>
\endxy}\Ea 
\right)=0
, \ \ 
\delta_v 
\left( \Ba{c}\resizebox{8mm}{!}{  \xy
(0,6)*{}="1";
    (0,0.2)*{\bu}="L";
  (-4,-5)*{}="C";
   (+4,-5)*{}="D";
\ar @{->} "D";"L" <0pt>
\ar @{->} "C";"L" <0pt>
\ar @{<-} "1";"L" <0pt>
 \endxy}
 \Ea\right)=0\ \ , \ \
\delta_v 
\left( \Ba{c}\resizebox{8mm}{!}{  \xy
(0,-6)*{}="1";
    (0,-0.2)*{\bu}="L";
  (-4,5)*{}="C";
   (+4,5)*{}="D";
\ar @{<-} "D";"L" <0pt>
\ar @{<-} "C";"L" <0pt>
\ar @{->} "1";"L" <0pt>
 \endxy}
 \Ea\right)=
   \Ba{c}\resizebox{11.5mm}{!}{  \xy
(0,6)*{}="1";
(8,0.2)*{}="2";
    (0,0.2)*{\bu}="L";
  (-4,-5)*{}="C";
   (4,-5)*{\bu}="D";
\ar @{->} "D";"L" <0pt>
\ar @{->} "C";"L" <0pt>
\ar @{<-} "1";"L" <0pt>
\ar @{<-} "2";"D" <0pt>
 \endxy}
 \Ea
 +
   \Ba{c}\resizebox{11.5mm}{!}{  \xy
(0,6)*{}="1";
(-8,0.2)*{}="2";
    (0,0.2)*{\bu}="L";
  (4,-5)*{}="C";
   (-4,-5)*{\bu}="D";
\ar @{->} "D";"L" <0pt>
\ar @{->} "C";"L" <0pt>
\ar @{<-} "1";"L" <0pt>
\ar @{<-} "2";"D" <0pt>
 \endxy}
 \Ea
\Eeq
\Beq\label{2: d on 4val}
\delta_v\left(  \Ba{c}\resizebox{17.5mm}{!}{  \xy
(0,-7)*+{_d}="D";
(0,+7)*+{_b}="U";
    (0,0)*{\bu}="C";
  (7,0)*+{_c}="R";
   (-7,0)*+{_a}="L";
\ar @{->} "C";"U" <0pt>
\ar @{->} "C";"D" <0pt>
\ar @{<-} "C";"L" <0pt>
\ar @{<-} "C";"R" <0pt>
 \endxy}
 \Ea\right)
 =
  \Ba{c}\resizebox{18.5mm}{!}{  \xy
(0,-7)*+{_d}="D";
(0,5)*{\bu}="U";
(0,+12)*+{_b}="0";
    (0,0)*{\bu}="C";
  (7,0)*+{_c}="R";
   (-7,5)*+{_a}="L";
\ar @{->} "U";"0" <0pt>  
\ar @{->} "C";"U" <0pt>
\ar @{->} "C";"D" <0pt>
\ar @{<-} "U";"L" <0pt>
\ar @{<-} "C";"R" <0pt>
 \endxy}
 \Ea 
 +
  \Ba{c}\resizebox{18.5mm}{!}{  \xy
(0,-7)*+{_b}="D";
(0,5)*{\bu}="U";
(0,+12)*+{_d}="0";
    (0,0)*{\bu}="C";
  (7,0)*+{_a}="R";
   (-7,5)*+{_c}="L";
\ar @{->} "U";"0" <0pt>  
\ar @{->} "C";"U" <0pt>
\ar @{->} "C";"D" <0pt>
\ar @{<-} "U";"L" <0pt>
\ar @{<-} "C";"R" <0pt>
 \endxy}
 \Ea
  +
  \Ba{c}\resizebox{18.5mm}{!}{  \xy
(0,-7)*+{_d}="D";
(0,5)*{\bu}="U";
(0,+12)*+{_b}="0";
    (0,0)*{\bu}="C";
  (7,5)*+{_c}="R";
   (-7,0)*+{_a}="L";
\ar @{->} "U";"0" <0pt>  
\ar @{->} "C";"U" <0pt>
\ar @{->} "C";"D" <0pt>
\ar @{<-} "U";"R" <0pt>
\ar @{<-} "C";"L" <0pt>
 \endxy}
 \Ea 
   +
  \Ba{c}\resizebox{18.5mm}{!}{  \xy
(0,-7)*+{_b}="D";
(0,5)*{\bu}="U";
(0,+12)*+{_d}="0";
    (0,0)*{\bu}="C";
  (7,5)*+{_a}="R";
   (-7,0)*+{_c}="L";
\ar @{->} "U";"0" <0pt>  
\ar @{->} "C";"U" <0pt>
\ar @{->} "C";"D" <0pt>
\ar @{<-} "U";"R" <0pt>
\ar @{<-} "C";"L" <0pt>
 \endxy}
 \Ea 
\Eeq

\subsubsection{\bf Passing vertices} A $(1,1)$-corolla $v=\Ba{c}\resizebox{3mm}{!}{ \xy
 (0,5)*{}="a",
(0,0)*{\bullet}="b",
(0,-5)*{}="c",
\ar @{<-} "a";"b" <0pt>
\ar @{<-} "b";"c" <0pt>
\endxy}\Ea$ is called a passing vertex. As the differential $\delta_v$ acts on such a vertex just by creating a second passing vertex,
$$
\delta_v \left(\Ba{c}\resizebox{3mm}{!}{ \xy
 (0,5)*{}="a",
(0,0)*{\bullet}="b",
(0,-5)*{}="c",
\ar @{<-} "a";"b" <0pt>
\ar @{<-} "b";"c" <0pt>
\endxy}\Ea\right)=\Ba{c}\resizebox{3mm}{!}{ \xy
 (0,5)*{\bu}="a",
(0,0)*{\bullet}="b",
(0,-5)*{}="c",
(0,10)*{}="u",
\ar @{<-} "a";"b" <0pt>
\ar @{<-} "b";"c" <0pt>
\ar @{<-} "u";"a" <0pt>
\endxy}\Ea
$$
it is easy to see that the differential ideal in $(\PCYf_d,\delta)$ generated by hairy ribbon graphs containing at least one passing vertex, is acyclic. Hence we can and will assume from now on that {\it the dg properad $(\PCYf_d,\delta)$ is spanned by graphs with no passing vertices}.

\subsection{A small model for $\PCYf_d$} Consider a non-differential ideal $J$ in the properad $\PCYf_d$ generated by (i) all ribbon $(m,0)$-corollas with $m\geq 3$, and (ii) all  ribbon $(m,n)$-corollas with $m+n\geq 4$ except the one given in (\ref{1: DLie corolla}).
Let $\langle J, \delta J\rangle$ be the differential closure of $J$ in the dg properad $\PCYf_d$. The quotient dg properad
$$
\Delta\PCYf_d:= \PCYf_d/\langle J, \delta J\rangle
$$
is generated by four ribbon corollas,
\Beq\label{1: types of corollas in DeltaPCY}
\Ba{c}\resizebox{8mm}{!}{\xy
 (0,0)*{}="a",
(4,-3)*{\bullet}="b",
(8,0)*{}="c",
\ar @{<-} "a";"b" <0pt>
\ar @{->} "b";"c" <0pt>
\endxy} \Ea, 
 \ \ \
 \Ba{c}\resizebox{8mm}{!}{  \xy
(0,-6)*{}="1";
    (0,-0.2)*{\bu}="L";
  (-4,5)*{}="C";
   (4,5)*{}="D";
\ar @{<-} "D";"L" <0pt>
\ar @{<-} "C";"L" <0pt>
\ar @{->} "1";"L" <0pt>
 \endxy}
 \Ea, 
 \ \ \ 
 \Ba{c}\resizebox{8mm}{!}{  \xy
(0,6)*{}="1";
    (0,0.2)*{\bu}="L";
  (-4,-5)*{}="C";
   (+4,-5)*{}="D";
\ar @{->} "D";"L" <0pt>
\ar @{->} "C";"L" <0pt>
\ar @{<-} "1";"L" <0pt>
 \endxy}
 \Ea,
 \ \ \  
  \Ba{c}\resizebox{13mm}{!}{  \xy
(0,-5)*{}="U";
(0,+5)*{}="D";
    (0,0)*{\bu}="C";
  (5,0)*{}="R";
   (-5,0)*{}="L";
\ar @{->} "C";"U" <0pt>
\ar @{->} "C";"D" <0pt>
\ar @{<-} "C";"L" <0pt>
\ar @{<-} "C";"R" <0pt>
 \endxy}
 \Ea\ ,
\Eeq
modulo the following quadratic relations for any internal edge in a ribbon graph from $\Delta\PCYf_d$,
\Beq\label{2: relations I in DeltaPCYf}
 \Ba{c}\resizebox{13.5mm}{!}{  \xy
    (-9,-3)*{\bbu}="L";
 (-14,3.5)*{\bbu}="B";
 (-20,12)*+{_1}="b1";
 (-8,12)*+{_2}="b2";
  (-3,4)*+{_3}="C";
\ar @{<-} "C";"L" <0pt>
\ar @{<-} "B";"L" <0pt>
\ar @{->} "B";"b1" <0pt>
\ar @{->} "B";"b2" <0pt>
 \endxy}
 \Ea
 +
  \Ba{c}\resizebox{13.5mm}{!}{  \xy
    (-9,-3)*{\bbu}="L";
 (-14,3.5)*{\bbu}="B";
 (-20,12)*+{_3}="b1";
 (-8,12)*+{_1}="b2";
  (-3,4)*+{_2}="C";
\ar @{<-} "C";"L" <0pt>
\ar @{<-} "B";"L" <0pt>
\ar @{->} "B";"b1" <0pt>
\ar @{->} "B";"b2" <0pt>
 \endxy}
 \Ea
 +
   \Ba{c}\resizebox{13.5mm}{!}{  \xy
    (-9,-3)*{\bbu}="L";
 (-14,3.5)*{\bbu}="B";
 (-20,12)*+{_2}="b1";
 (-8,12)*+{_3}="b2";
  (-3,4)*+{_1}="C";
\ar @{<-} "C";"L" <0pt>
\ar @{<-} "B";"L" <0pt>
\ar @{->} "B";"b1" <0pt>
\ar @{->} "B";"b2" <0pt>
 \endxy}
 \Ea
 \hspace{-1mm}
 =0\ , 
 \ \
   \Ba{c}\resizebox{9.5mm}{!}{  \xy
(-5,12)*+{_1}="1";
 (5,12)*+{_2}="2";
    (0,3.5)*{\bbu}="A";
 (0,-3.5)*{\bbu}="B";
 (-5,-12)*+{_4}="b1";
 (5,-12)*+{_3}="b2";
\ar @{->} "A";"1" <0pt>
\ar @{->} "A";"2" <0pt>
\ar @{->} "B";"A" <0pt>
\ar @{<-} "B";"b1" <0pt>
\ar @{<-} "B";"b2" <0pt>
 \endxy}
 \Ea
 +
 \Ba{c}\resizebox{13.5mm}{!}{  \xy
(-9,12)*+{_2}="1";
    (-9,+3)*{\bbu}="L";
 (-14,-3.5)*{\bbu}="B";
 (-19,5)*+{_1}="b1";
 (-14,-12)*+{_4}="b2";
  (-3,-5)*+{_3}="C";
\ar @{->} "C";"L" <0pt>
\ar @{->} "B";"L" <0pt>
\ar @{->} "B";"b1" <0pt>
\ar @{<-} "B";"b2" <0pt>
\ar @{<-} "1";"L" <0pt>
 \endxy}
 \Ea
 +
 \Ba{c}\resizebox{13.5mm}{!}{  \xy
(-9,-10)*+{_3}="1";
    (-9,-3)*{\bbu}="L";
 (-14,3.5)*{\bbu}="B";
 (-19,-5)*+{_4}="b1";
 (-14,12)*+{_1}="b2";
  (-3,5)*+{_2}="C";
\ar @{<-} "C";"L" <0pt>
\ar @{<-} "B";"L" <0pt>
\ar @{<-} "B";"b1" <0pt>
\ar @{->} "B";"b2" <0pt>
\ar @{->} "1";"L" <0pt>
 \endxy}
 \Ea
=0
\Eeq
\Beq\label{2: twisted ass-type relations}
 \Ba{c}\resizebox{13.5mm}{!}{  \xy
(-9,10)*{^0}="1";
    (-9,+3)*{\bbu}="L";
 (-14,-3.5)*{\bbu}="B";
 (-20,-12)*+{_1}="b1";
 (-8,-12)*+{_2}="b2";
  (-3,-4)*{_3}="C";
\ar @{->} "C";"L" <0pt>
\ar @{->} "B";"L" <0pt>
\ar @{<-} "B";"b1" <0pt>
\ar @{<-} "B";"b2" <0pt>
\ar @{<-} "1";"L" <0pt>
 \endxy}
 \Ea
 +
  \Ba{c}\resizebox{13.5mm}{!}{  \xy
(9,10)*{^0}="1";
    (9,+3)*{\bbu}="L";
 (14,-3.5)*{\bbu}="B";
 (20,-12)*+{_3}="b1";
 (8,-12)*+{2}="b2";
  (3,-4)*{1}="C";
\ar @{->} "C";"L" <0pt>
\ar @{->} "B";"L" <0pt>
\ar @{<-} "B";"b1" <0pt>
\ar @{<-} "B";"b2" <0pt>
\ar @{<-} "1";"L" <0pt>
 \endxy}
 \Ea
 =0, 
 \ \ \
 \Ba{c}\resizebox{13.5mm}{!}{  \xy
(-9,-10)*{_0}="1";
    (-9,-3)*{\bbu}="L";
 (-14,3.5)*{\bbu}="B";
 (-20,12)*+{_1}="b1";
 (-8,12)*+{_2}="b2";
  (-3,4)*{_3}="C";
\ar @{<-} "C";"L" <0pt>
\ar @{<-} "B";"L" <0pt>
\ar @{->} "B";"b1" <0pt>
\ar @{->} "B";"b2" <0pt>
\ar @{->} "1";"L" <0pt>
 \endxy}
 \Ea
 +
  \Ba{c}\resizebox{13.5mm}{!}{  \xy
(9,-10)*{_0}="1";
    (9,-3)*{\bbu}="L";
 (14,3.5)*{\bbu}="B";
 (20,12)*+{_3}="b1";
 (8,12)*+{_2}="b2";
  (3,4)*{_1}="C";
\ar @{<-} "C";"L" <0pt>
\ar @{<-} "B";"L" <0pt>
\ar @{->} "B";"b1" <0pt>
\ar @{->} "B";"b2" <0pt>
\ar @{->} "1";"L" <0pt>
 \endxy}
 \Ea
 -
    \Ba{c}\resizebox{15mm}{!}{  \xy
 (0,14)*{_2}="1";
(0,-7)*{^0}="U";
(0,7)*{\bu}="D";
    (0,0)*{\bu}="C";
  (7,0)*{_3}="R";
   (-7,0)*{_1}="L";
\ar @{<-} "C";"U" <0pt>
\ar @{<-} "C";"D" <0pt>
\ar @{->} "C";"L" <0pt>
\ar @{->} "C";"R" <0pt>
\ar @{->} "D";"1" <0pt>
 \endxy}
 \Ea
 =0,
\Eeq

\Beq\label{2: 3+4 relations}
  \Ba{c}\resizebox{17mm}{!}{  \xy
(-3,13)*+{_1}="UL";
(3,13)*+{_2}="UR";
(0,7)*{\bu}="U";
(0,-7)*+{_4}="D";
    (0,0)*{\bu}="C";
  (7,0)*+{_3}="R";
   (-7,0)*+{_0}="L";
\ar @{->} "U";"UL" <0pt>   
\ar @{->} "U";"UR" <0pt>  
\ar @{->} "C";"U" <0pt>
\ar @{->} "C";"D" <0pt>
\ar @{<-} "C";"L" <0pt>
\ar @{<-} "C";"R" <0pt>
 \endxy}
 \Ea
 \hspace{-1mm}
 -
 \hspace{-2mm}
   \Ba{c}\resizebox{17mm}{!}{  \xy
(-3,13)*+{_0}="UL";
(3,13)*+{_1}="UR";
(0,7)*{\bu}="U";
(0,-7)*+{_3}="D";
    (0,0)*{\bu}="C";
  (7,0)*+{_2}="R";
   (-7,0)*+{_4}="L";
\ar @{<-} "U";"UL" <0pt>   
\ar @{->} "U";"UR" <0pt>  
\ar @{<-} "C";"U" <0pt>
\ar @{<-} "C";"D" <0pt>
\ar @{->} "C";"L" <0pt>
\ar @{->} "C";"R" <0pt>
 \endxy}
 \Ea
 \hspace{-1mm}
 -
 \hspace{-2mm}
    \Ba{c}\resizebox{17mm}{!}{  \xy
(-3,13)*+{_2}="UL";
(3,13)*+{_3}="UR";
(0,7)*{\bu}="U";
(0,-7)*+{_0}="D";
    (0,0)*{\bu}="C";
  (7,0)*+{_4}="R";
   (-7,0)*+{_1}="L";
\ar @{->} "U";"UL" <0pt>   
\ar @{<-} "U";"UR" <0pt>  
\ar @{<-} "C";"U" <0pt>
\ar @{<-} "C";"D" <0pt>
\ar @{->} "C";"L" <0pt>
\ar @{->} "C";"R" <0pt>
 \endxy}
 \Ea
 \hspace{-2mm}
 =0, 
   \Ba{c}\resizebox{17mm}{!}{  \xy
(-3,13)*+{_1}="UL";
(3,13)*+{_2}="UR";
(0,7)*{\bu}="U";
(0,-7)*+{_4}="D";
    (0,0)*{\bu}="C";
  (7,0)*+{_3}="R";
   (-7,0)*+{_0}="L";
\ar @{<-} "U";"UL" <0pt>   
\ar @{<-} "U";"UR" <0pt>  
\ar @{<-} "C";"U" <0pt>
\ar @{<-} "C";"D" <0pt>
\ar @{->} "C";"L" <0pt>
\ar @{->} "C";"R" <0pt>
 \endxy}
 \Ea
 \hspace{-1mm}
 -
 \hspace{-2mm}
   \Ba{c}\resizebox{17mm}{!}{  \xy
(-3,13)*+{_0}="UL";
(3,13)*+{_1}="UR";
(0,7)*{\bu}="U";
(0,-7)*+{_3}="D";
    (0,0)*{\bu}="C";
  (7,0)*+{_2}="R";
   (-7,0)*+{_4}="L";
\ar @{->} "U";"UL" <0pt>   
\ar @{<-} "U";"UR" <0pt>  
\ar @{->} "C";"U" <0pt>
\ar @{->} "C";"D" <0pt>
\ar @{<-} "C";"L" <0pt>
\ar @{<-} "C";"R" <0pt>
 \endxy}
 \hspace{0mm}
-
\hspace{-2mm}
     \Ba{c}\resizebox{17mm}{!}{  \xy
(-3,13)*+{^2}="UL";
(3,13)*+{_3}="UR";
(0,7)*{\bu}="U";
(0,-7)*+{_0}="D";
    (0,0)*{\bu}="C";
  (7,0)*+{_4}="R";
   (-7,0)*+{_1}="L";
\ar @{<-} "U";"UL" <0pt>   
\ar @{->} "U";"UR" <0pt>  
\ar @{->} "C";"U" <0pt>
\ar @{->} "C";"D" <0pt>
\ar @{<-} "C";"L" <0pt>
\ar @{<-} "C";"R" <0pt>
 \endxy}
 \Ea
 \Ea
 \hspace{-2mm}
 =0,
\Eeq
\Beq\label{2: double Lie relations}
  \Ba{c}\resizebox{26mm}{!}{  \xy
(0,-7)*+{_5}="D";
(0,+7)*+{_1}="U";
 (-7,0)*+{_0}="L";
(0,0)*{\bu}="C";
(7,0)*{\bu}="CR";
  (7,-7)*+{_4}="RD";
(7,+7)*+{_2}="RU";
 (14,0)*+{_3}="RR";
\ar @{->} "C";"U" <0pt>
\ar @{->} "C";"D" <0pt>
\ar @{<-} "C";"L" <0pt>
\ar @{<-} "C";"CR" <0pt>
\ar @{<-} "CR";"RU" <0pt>
\ar @{<-} "CR";"RD" <0pt>
\ar @{->} "CR";"RR" <0pt>
 \endxy}
 \Ea
 +
  \Ba{c}\resizebox{26mm}{!}{  \xy
(0,-7)*+{_3}="D";
(0,+7)*+{_5}="U";
 (-7,0)*+{_4}="L";
(0,0)*{\bu}="C";
(7,0)*{\bu}="CR";
  (7,-7)*+{_2}="RD";
(7,+7)*+{_0}="RU";
 (14,0)*+{_1}="RR";
\ar @{->} "C";"U" <0pt>
\ar @{->} "C";"D" <0pt>
\ar @{<-} "C";"L" <0pt>
\ar @{<-} "C";"CR" <0pt>
\ar @{<-} "CR";"RU" <0pt>
\ar @{<-} "CR";"RD" <0pt>
\ar @{->} "CR";"RR" <0pt>
 \endxy}
 \Ea
  +
  \Ba{c}\resizebox{26mm}{!}{  \xy
(0,-7)*+{_1}="D";
(0,+7)*+{_3}="U";
 (-7,0)*+{_2}="L";
(0,0)*{\bu}="C";
(7,0)*{\bu}="CR";
  (7,-7)*+{_0}="RD";
(7,+7)*+{_4}="RU";
 (14,0)*+{_5}="RR";
\ar @{->} "C";"U" <0pt>
\ar @{->} "C";"D" <0pt>
\ar @{<-} "C";"L" <0pt>
\ar @{<-} "C";"CR" <0pt>
\ar @{<-} "CR";"RU" <0pt>
\ar @{<-} "CR";"RD" <0pt>
\ar @{->} "CR";"RR" <0pt>
 \endxy}
 \Ea=0.
\Eeq
The labelled arrows in the above pictures can stand either for hairs or half-edges.  
The sign rule is determined by our definition of the orientations of hairy 
ribbon graphs given in (\ref{2: or for d even}) and (\ref{2: or for d odd}), and by the agreement 
all hairs and half edges are ordered in accordance with their labels. 
For $d$  odd (when the internal edges have even degree, while vertices have odd degree) we always assume in our pictures that vertices are ordered from the top to the bottom (as in  (\ref{2: relations I in DeltaPCYf}) and 
(\ref{2: 3+4 relations}), and, if vertices are on the same level  (as in (\ref{2: double Lie relations})), they are ordered from the left to the right.

\subsection{Proof of Main Theorem {\ref{1: Theorem on PCYf}}}
The main idea of the proof is to use appropriate filtrations of the complexes $\PCYf_d(m,n)$ and $\Delta\PCYf_d(m,n)$ such that on the first two pages of the associated spectral sequences the problem reduces essentially to the study of the tensor products of complexes whose cohomologies are already known (they are listed in the next subsection). This scenario is close  to the one used in \cite{Me2} in the study of a low-valence complex of ribbon quivers which computes the cohomology of the moduli spaces $\cM_{g,m}$, but there are some deviations. Hence we present the proof in full detail.

\subsubsection{\sf Several auxiliary complexes whose cohomology is well-known}\label{2: subsec on aux complexes}

  {\bf (I)} The dg free operad  $\cA ss_\infty=\{\cA ss_\infty(n)\}_{n\geq 2}$ is generated by degree $2-n$ planar corollas
$$
  \underbrace{\Ba{c}\resizebox{20mm}{!}
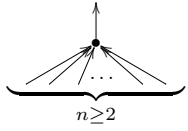
{\xy
(1,-6)*{\ldots},
(-13,-7)*{},
(-8,-7)*{},
(-3,-7)*{},
(7,-7)*{},
(13,-7)*{},
 (0,0)*{\bu}="a",
(0,7)*{}="0",
(-12,-7)*{}="b_1",
(-8,-7)*{}="b_2",
(-3,-7)*{}="b_3",
(8,-7)*{}="b_4",
(12,-7)*{}="b_5",
\ar @{->} "a";"0" <0pt>
\ar @{<-} "a";"b_2" <0pt>
\ar @{<-} "a";"b_3" <0pt>
\ar @{<-} "a";"b_1" <0pt>
\ar @{<-} "a";"b_4" <0pt>
\ar @{<-} "a";"b_5" <0pt>
\endxy}\Ea}_{n\geq 2}
$$
 with $n\geq 2$ incoming totally ordered legs and one outgoing leg. The differential is given on the generators by the standard sum $\sum_v\delta_v$ over all vertices, where $\delta_v$ means 
 substituting  the graph $\Ba{c}\xy
 (0,0)*{\bu}="a",
(0,5)*{\bu}="b",
\ar @{->} "a";"b" <0pt>
\endxy\Ea$ into the vertex, and reattaching the edges among the two newly created vertices in such a way that the total order of incoming edges is preserved. Thus one has
\Beq\label{2: d in Ass_infty}
\delta_v
  \underbrace{\Ba{c}\resizebox{20mm}{!}{\xy
(1,-6)*{\ldots},
(-13,-7)*{},
(-8,-7)*{},
(-3,-7)*{},
(7,-7)*{},
(13,-7)*{},
 (0,0)*{\bu}="a",
(0,7)*{}="0",
(-12,-7)*{}="b_1",
(-8,-7)*{}="b_2",
(-3,-7)*{}="b_3",
(8,-7)*{}="b_4",
(12,-7)*{}="b_5",
\ar @{->} "a";"0" <0pt>
\ar @{<-} "a";"b_2" <0pt>
\ar @{<-} "a";"b_3" <0pt>
\ar @{<-} "a";"b_1" <0pt>
\ar @{<-} "a";"b_4" <0pt>
\ar @{<-} "a";"b_5" <0pt>
\endxy}\Ea}_n
=
\sum_{A\subsetneq [n]} 
\Ba{c}\resizebox{22mm}{!}{\xy
(1.6,-7)*{...},
(-13,-7)*{},
(-8,-7)*{},
(-3,-7)*{},
(7,-7)*{},
(13,-7)*{},
 (0,0)*{\bu}="a",
(0,7)*{}="0",
(-12,-7)*{}="b_1",
(-9,-7)*{}="b_2",
(-6,-7)*{...},
(-3,-7)*{\bu}="b_3",
(-9,-14)*{}="c1",
(-5,-14)*{}="c2",
(-1,-14)*{...},
(3,-14)*{}="c3",
(7,-7)*{}="b_4",
(11,-7)*{}="b_5",
(-3,-17)*{\underbrace{\ \ \ \ \ \ \ \ \ \ \ \ }_A},
\ar @{->} "a";"0" <0pt>
\ar @{<-} "a";"b_2" <0pt>
\ar @{<-} "a";"b_3" <0pt>
\ar @{<-} "a";"b_1" <0pt>
\ar @{<-} "a";"b_4" <0pt>
\ar @{<-} "a";"b_5" <0pt>
\ar @{->} "c1";"b_3" <0pt>
\ar @{->} "c2";"b_3" <0pt>
\ar @{->} "c3";"b_3" <0pt>
\endxy}\Ea
\Eeq
where the summation runs over proper connected  subsets $A$ (of cardinality $\geq 2$) of the totally ordered set $[n]$.
Its cohomology $H^\bu(\cA ss_\infty)$ was proven in \cite{GK} to be the operad $\cA ss$ generated by the planar corolla $
  \Ba{c}\resizebox{8mm}{!}{  \xy
(0,6)*{}="1";
    (0,0.2)*{\bu}="L";
  (-4,-5)*{}="C";
   (+4,-5)*{}="D";
\ar @{->} "D";"L" <0pt>
\ar @{->} "C";"L" <0pt>
\ar @{<-} "1";"L" <0pt>
 \endxy}
 \Ea
$ (whose two incoming legs are totally ordered from the left to the right) modulo the associativity relation
 \Beq\label{3: Ass operad relation}
  \Ba{c}\resizebox{13.5mm}{!}{  \xy
(-9,10)*{^0}="1";
    (-9,+3)*{\bbu}="L";
 (-14,-3.5)*{\bbu}="B";
 (-20,-12)*+{_1}="b1";
 (-8,-12)*+{_2}="b2";
  (-3,-4)*{_3}="C";
\ar @{->} "C";"L" <0pt>
\ar @{->} "B";"L" <0pt>
\ar @{<-} "B";"b1" <0pt>
\ar @{<-} "B";"b2" <0pt>
\ar @{<-} "1";"L" <0pt>
 \endxy}
 \Ea
 +
  \Ba{c}\resizebox{13.5mm}{!}{  \xy
(9,10)*{^0}="1";
    (9,+3)*{\bbu}="L";
 (14,-3.5)*{\bbu}="B";
 (20,-12)*+{_3}="b1";
 (8,-12)*+{2}="b2";
  (3,-4)*{1}="C";
\ar @{->} "C";"L" <0pt>
\ar @{->} "B";"L" <0pt>
\ar @{<-} "B";"b1" <0pt>
\ar @{<-} "B";"b2" <0pt>
\ar @{<-} "1";"L" <0pt>
 \endxy}
 \Ea=0.
 \Eeq
 Denote by $\cA ss_\infty^{op}$ and $\cA ss^{op}$ the versions of the above two operads in which the directions of all arrows are reversed.

 \sip

{\bf (II)}
   Let $\cA ss_\infty^{\wedge, +}=\{\cA ss_\infty^{\wedge, +}(n)\}_{n\geq 1}$  be a dg free right module over the dg operad $\cA ss_\infty$ generated by $(0,n)$-corollas 
 $$
  \underbrace{\Ba{c}\resizebox{20mm}{!}
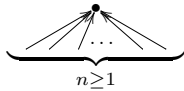
{\xy
(1,-6)*{\ldots},
(-13,-7)*{},
(-8,-7)*{},
(-3,-7)*{},
(7,-7)*{},
(13,-7)*{},
 (0,0)*{\bu}="a",
(0,7)*{}="0",
(-12,-7)*{}="b_1",
(-8,-7)*{}="b_2",
(-3,-7)*{}="b_3",
(8,-7)*{}="b_4",
(12,-7)*{}="b_5",
%
\ar @{<-} "a";"b_2" <0pt>
\ar @{<-} "a";"b_3" <0pt>
\ar @{<-} "a";"b_1" <0pt>
\ar @{<-} "a";"b_4" <0pt>
\ar @{<-} "a";"b_5" <0pt>
\endxy}\Ea}_{n\geq 1}
$$
whose incoming legs are {\it totally ordered}\, from the left to the right; they are assigned the cohomological degree $d-n$. The differential $\delta^+$ in  $\cA ss_\infty^{\wedge, +}$ is given by the standard sum over vertices 
$\delta^+=\sum_v\delta^+_v $, where $\delta_v^+$ acts on $(1,n)$ vertices as in  (\ref{2: d in Ass_infty}), while on $(0,n)$-vertices $v$ as follows, 
\Beq\label{2: d in Ass_wedge_plus_infty}
\delta^+ \hspace{-2mm}  \Ba{c}\resizebox{2.5mm}{!}{  \xy
    (0,2)*{\bu}="L";
  (0,-4)*{}="C";
\ar @{->} "C";"L" <0pt>
 \endxy}
 \Ea=0, \ \ \
 \delta^+ \hspace{-3mm}   \Ba{c}\resizebox{8mm}{!}{  \xy
(0,6)*{}="1";
    (0,0.2)*{\bu}="L";
  (-4,-5)*{}="C";
   (+4,-5)*{}="D";
\ar @{->} "D";"L" <0pt>
\ar @{->} "C";"L" <0pt>
%
 \endxy}
 \Ea
 =  \Ba{c}\resizebox{8mm}{!}{  \xy
(0,6)*{\bu}="1";
    (0,0.2)*{\bu}="L";
  (-4,-5)*{}="C";
   (+4,-5)*{}="D";
\ar @{->} "D";"L" <0pt>
\ar @{->} "C";"L" <0pt>
\ar @{<-} "1";"L" <0pt>
 \endxy}
 \Ea, \ \ \ \ \
 \delta^+\hspace{-3mm}
  \underbrace{\Ba{c}\resizebox{20mm}{!}{\xy
(1,-6)*{\ldots},
(-13,-7)*{},
(-8,-7)*{},
(-3,-7)*{},
(7,-7)*{},
(13,-7)*{},
 (0,0)*{\bu}="a",
(0,7)*{}="0",
(-12,-7)*{}="b_1",
(-8,-7)*{}="b_2",
(-3,-7)*{}="b_3",
(8,-7)*{}="b_4",
(12,-7)*{}="b_5",
%
\ar @{<-} "a";"b_2" <0pt>
\ar @{<-} "a";"b_3" <0pt>
\ar @{<-} "a";"b_1" <0pt>
\ar @{<-} "a";"b_4" <0pt>
\ar @{<-} "a";"b_5" <0pt>
\endxy}\Ea}_{n\geq 3}
=
\Ba{c}\resizebox{20mm}{!}{\xy
(1,-6)*{\ldots},
(-13,-7)*{},
(-8,-7)*{},
(-3,-7)*{},
(7,-7)*{},
(13,-7)*{},
 (0,0)*{\bu}="a",
(0,7)*{\bu}="0",
(-12,-7)*{}="b_1",
(-8,-7)*{}="b_2",
(-3,-7)*{}="b_3",
(8,-7)*{}="b_4",
(12,-7)*{}="b_5",
\ar @{->} "a";"0" <0pt>
\ar @{<-} "a";"b_2" <0pt>
\ar @{<-} "a";"b_3" <0pt>
\ar @{<-} "a";"b_1" <0pt>
\ar @{<-} "a";"b_4" <0pt>
\ar @{<-} "a";"b_5" <0pt>
\endxy}\Ea\hspace{-3mm}
+
\sum_{A\subsetneq [n]\atop \# A\geq 2} 
\Ba{c}\resizebox{22mm}{!}{\xy
(1.6,-7)*{...},
(-13,-7)*{},
(-8,-7)*{},
(-3,-7)*{},
(7,-7)*{},
(13,-7)*{},
 (0,0)*{\bu}="a",
(0,7)*{}="0",
(-12,-7)*{}="b_1",
(-9,-7)*{}="b_2",
(-6,-7)*{...},
(-3,-7)*{\bu}="b_3",
(-9,-14)*{}="c1",
(-5,-14)*{}="c2",
(-1,-14)*{...},
(3,-14)*{}="c3",
(7,-7)*{}="b_4",
(11,-7)*{}="b_5",
(-3,-17)*{\underbrace{\ \ \ \ \ \ \ \ \ \ \ \ }_A},
%
\ar @{<-} "a";"b_2" <0pt>
\ar @{<-} "a";"b_3" <0pt>
\ar @{<-} "a";"b_1" <0pt>
\ar @{<-} "a";"b_4" <0pt>
\ar @{<-} "a";"b_5" <0pt>
\ar @{->} "c1";"b_3" <0pt>
\ar @{->} "c2";"b_3" <0pt>
\ar @{->} "c3";"b_3" <0pt>
\endxy}\Ea.
\Eeq
Using a filtration of each complex $\cA ss_\infty^{\wedge,+}(n)$, $n\geq 1$, by the number of non-univalent vertices, one concludes immediately from the associated graded complex that
\Beq\label{3: H(Ass_infty_wedge+plus)}
H^\bu(\cA ss_\infty^{\wedge,+}(1))=\K[d-1], \ \ \ \ H^\bu(\cA ss_\infty^{\wedge,+}(n))=0\ \ 
\text{for 
$n\geq 2$},
\Eeq
i.e.\ the cohomology of $\cA ss_\infty^{\wedge,+}$ is 1-dimensional and is spanned by \hspace{-2mm}
$ \Ba{c}\resizebox{2.5mm}{!}{  \xy
    (0,2)*{\bu}="L";
  (0,-4)*{}="C";
\ar @{->} "C";"L" <0pt>
 \endxy}
 \Ea$.
Let us denote by $\cA ss_\infty^{\vee,+}$ a copy of $\cA ss_\infty^{\wedge,+}$ 
in which directions of all arrows in the generators are reversed. 

\sip

{\bf (III)}  Let $\cA ss_\infty^{\vee, cyc}=\{\cA ss_\infty^{\vee,cyc}(n)\}_{n\geq 2}$ be a dg free right module over the dg operad $\cA ss_\infty$ generated by degree $(2-d)m +d$ ribbon corollas of the form
 \Beq\label{2: sources}
\Ba{c}\resizebox{20mm}{!}{\xy
(2,-6)*{\ldots},
(-13,+7)*{},
(-8,-1)*{},
(-3,-7)*{},
(7,-7)*{},
(13,-7)*{},
 (0,0)*{\bu}="a",
(0,7)*{}="0",
(-7,+7)*+{_m}="b_1",
(-9,-3)*+{_1}="b_2",
(-3,-9)*+{_2}="b_3",
(8,-7)*+{}="b_4",
(5,9)*+{_{m-1}}="b_5",
%
\ar @{->} "a";"b_2" <0pt>
\ar @{->} "a";"b_3" <0pt>
\ar @{->} "a";"b_1" <0pt>
\ar @{->} "a";"b_4" <0pt>
\ar @{->} "a";"b_5" <0pt>
\endxy}\Ea, \ \ \ {k\geq 2}.
\Eeq
Vertices of ribbon graphs having this form are called {\it sources}.
 
\sip

Let $((n))$ denote the finite set   $\{1,2,...,n\}$ equipped with the cyclic ordering  $1<2<...<n<1$. The differential in $\cA ss_\infty^{\vee, cyc}$ is defined by
\Beq\label{2: d in Ass_wedge_infty}
\delta
\Ba{c}\resizebox{20mm}{!}{\xy
(2,-6)*{\ldots},
(-13,+7)*{},
(-8,-1)*{},
(-3,-7)*{},
(7,-7)*{},
(13,-7)*{},
 (0,0)*{\bu}="a",
(0,7)*{}="0",
(-7,+7)*+{_n}="b_1",
(-9,-3)*+{_1}="b_2",
(-3,-9)*+{_2}="b_3",
(8,-7)*+{}="b_4",
(6,7)*+{}="b_5",
%
\ar @{->} "a";"b_2" <0pt>
\ar @{->} "a";"b_3" <0pt>
\ar @{->} "a";"b_1" <0pt>
\ar @{->} "a";"b_4" <0pt>
\ar @{->} "a";"b_5" <0pt>
\endxy}\Ea
=
\sum_{A\subsetneq ((n))\atop \# A\geq 2} 
\Ba{c}\resizebox{20mm}{!}{\xy
(-13,-7)*{},
(-8,-7)*{},
(-3,-7)*{},
(7,-7)*{},
(13,-7)*{},
 (0,0)*{\bu}="a",
(0,7)*{}="0",
(-6,8)*{}="b_1",
(-9,-5)*{}="b_2",
(-3,-7)*{\bu}="b_3",
(-9,-14)*{}="c1",
(-5,-14)*{}="c2",
(-1,-14)*{...},
(3,-14)*{}="c3",
(7,-5)*{}="b_4",
(6,8)*{}="b_5",
(-3,-17)*{\underbrace{\ \ \ \ \ \ \ \ \ \ \ \ }_A},
%
\ar @{->} "a";"b_2" <0pt>
\ar @{->} "a";"b_3" <0pt>
\ar @{->} "a";"b_1" <0pt>
\ar @{->} "a";"b_4" <0pt>
\ar @{->} "a";"b_5" <0pt>
\ar @{<-} "c1";"b_3" <0pt>
\ar @{<-} "c2";"b_3" <0pt>
\ar @{<-} "c3";"b_3" <0pt>
\endxy}\Ea,
\Eeq
where the summation runs over proper connected subsets $A$ of $((n))$ of cardinality $\geq 2$.  For example
$$
\delta \Ba{c}\resizebox{16mm}{!}{ \xy
 (-7,2)*+{_1}="1",
 (0,0)*+{_2}="2",
 (7,2)*+{_3}="3",
(0,8)*{\bu}="c",
\ar @{->} "c";"1" <0pt>
\ar @{->} "c";"2" <0pt>
\ar @{->} "c";"3" <0pt>
\endxy}\Ea
=
  \Ba{c}\resizebox{13.5mm}{!}{  \xy
    (-9,+3)*{\bu}="L";
 (-14,-3.5)*{\bu}="B";
 (-20,-12)*+{_1}="b1";
 (-8,-12)*+{_2}="b2";
  (-3,-4)*{_3}="C";
\ar @{<-} "C";"L" <0pt>
\ar @{<-} "B";"L" <0pt>
\ar @{->} "B";"b1" <0pt>
\ar @{->} "B";"b2" <0pt>
 \endxy}
 \Ea
 +
  \Ba{c}\resizebox{13.5mm}{!}{  \xy
    (-9,+3)*{\bu}="L";
 (-14,-3.5)*{\bu}="B";
 (-20,-12)*+{_2}="b1";
 (-8,-12)*+{_3}="b2";
  (-3,-4)*{_1}="C";
\ar @{<-} "C";"L" <0pt>
\ar @{<-} "B";"L" <0pt>
\ar @{->} "B";"b1" <0pt>
\ar @{->} "B";"b2" <0pt>
 \endxy}
 \Ea
 +
  \Ba{c}\resizebox{13.5mm}{!}{  \xy
    (-9,+3)*{\bu}="L";
 (-14,-3.5)*{\bu}="B";
 (-20,-12)*+{_3}="b1";
 (-8,-12)*+{_1}="b2";
  (-3,-4)*{_2}="C";
\ar @{<-} "C";"L" <0pt>
\ar @{<-} "B";"L" <0pt>
\ar @{->} "B";"b1" <0pt>
\ar @{->} "B";"b2" <0pt>
 \endxy}
 \Ea.
 $$

It is not hard to show (see \S 3.4 in \cite{Me1} for a one-page proof) that the cohomology $H^\bu(\cA ss_\infty^{\vee,cyc})$ equals  a right module 
$\cA ss^{\vee,cyc}$ over the operad $\cA ss$ which, by definition, is generated by the ribbon corolla
 $
  \Ba{c}\resizebox{8mm}{!}{  \xy
    (0,0.2)*{\circ}="L";
  (-4,-5)*{}="C";
   (+4,-5)*{}="D";
\ar @{->} "D";"L" <0pt>
\ar @{->} "C";"L" <0pt>
 \endxy}
 \Ea
$
 modulo the relation given by the first formula in (\ref{2: relations I in DeltaPCYf}).

{\bf (IV)}   Consider a dg free properad $\mathcal{IB}^{path}_\infty$ of generated  ribbon $(m,n)$-corollas \cite{A}
\Beq\label{3: (m,n) corolla}
\Ba{c}\resizebox{11mm}{!}{ \xy
(0,8)*{\overbrace{\  \ \ \ \ \ \ \ \ \ }^{m}},
(0,-8)*{\underbrace{\  \ \ \ \ \ \ \ \ \ }_{n}},
(0,4.5)*+{...},
(0,-4.5)*+{...},
(0,0)*{\bu}="o",
(-5,5)*{}="1",
(-3,5)*{}="2",
(3,5)*{}="3",
(5,5)*{}="4",
(-3,-5)*{}="5",
(3,-5)*{}="6",
(5,-5)*{}="7",
(-5,-5)*{}="8",
\ar @{->} "o";"1" <0pt>
\ar @{->} "o";"2" <0pt>
\ar @{->} "o";"3" <0pt>
\ar @{->} "o";"4" <0pt>
\ar @{<-} "o";"5" <0pt>
\ar @{<-} "o";"6" <0pt>
\ar @{<-} "o";"7" <0pt>
\ar @{<-} "o";"8" <0pt>
\endxy}\Ea
\Eeq
 with $m\geq 1$, $n\geq 1$, $m+n\geq 3$. Note that the number of paths connecting incoming hairs to outgoing ones is equal to $mn$.
The differential $\delta_0$ in $\mathcal{IB}^{path}_\infty$ acts on such an $(m,n)$-corolla by substituting into its vertex 
$\bu$ the graph  $\resizebox{9mm}{!}{ \xy
 (0,1)*{\bu}="a",
(6,1)*{\bu}="b",
\ar @{->} "a";"b" <0pt>
\endxy}$ and reattaching the edges among the two newly created vertices in such a way that the cyclic structure is preserved and the total number of directed paths  connecting the incoming hairs to outgoing ones stays equal to $mn$. Thus
\Beq\label{2: d in IB_infty_path}
\delta_0
\Ba{c}\resizebox{11mm}{!}{ \xy
(0,8)*{\overbrace{\  \ \ \ \ \ \ \ \ \ }^{m\geq 2}},
(0,-8)*{\underbrace{\  \ \ \ \ \ \ \ \ \ }_{n\geq 2}},
(0,4.5)*+{...},
(0,-4.5)*+{...},
(0,0)*{\bu}="o",
(-5,5)*{}="1",
(-3,5)*{}="2",
(3,5)*{}="3",
(5,5)*{}="4",
(-3,-5)*{}="5",
(3,-5)*{}="6",
(5,-5)*{}="7",
(-5,-5)*{}="8",
\ar @{->} "o";"1" <0pt>
\ar @{->} "o";"2" <0pt>
\ar @{->} "o";"3" <0pt>
\ar @{->} "o";"4" <0pt>
\ar @{<-} "o";"5" <0pt>
\ar @{<-} "o";"6" <0pt>
\ar @{<-} "o";"7" <0pt>
\ar @{<-} "o";"8" <0pt>
\endxy}\Ea
=
\Ba{c}\resizebox{11mm}{!}{ \xy
(0,9.5)*+{...},
(0,-4.5)*+{...},
(0,5)*{\bu}="oo",
(0,0)*{\bu}="o",
(-5,10)*{}="1",
(-3,10)*{}="2",
(3,10)*{}="3",
(5,10)*{}="4",
(-3,-5)*{}="5",
(3,-5)*{}="6",
(5,-5)*{}="7",
(-5,-5)*{}="8",
\ar @{->} "oo";"1" <0pt>
\ar @{->} "oo";"2" <0pt>
\ar @{->} "oo";"3" <0pt>
\ar @{->} "oo";"4" <0pt>
\ar @{->} "o";"oo" <0pt>
\ar @{<-} "o";"5" <0pt>
\ar @{<-} "o";"6" <0pt>
\ar @{<-} "o";"7" <0pt>
\ar @{<-} "o";"8" <0pt>
\endxy}\Ea
+
\sum_{A\subsetneq [n]\atop \# A\geq 2} 
\Ba{c}\resizebox{22mm}{!}{\xy
(1.6,-7)*{...},
(-13,-7)*{},
(-8,-7)*{},
(-3,-7)*{},
(7,-7)*{},
(13,-7)*{},
 (0,0)*{\bu}="a",
 (0,6)*{...},
(-5,7)*{}="01",
(-3,7)*{}="02",
(3,7)*{}="03",
(5,7)*{}="04",
(-12,-7)*{}="b_1",
(-9,-7)*{}="b_2",
(-6,-7)*{...},
(-3,-7)*{\bu}="b_3",
(-9,-14)*{}="c1",
(-5,-14)*{}="c2",
(-1,-14)*{...},
(3,-14)*{}="c3",
(7,-7)*{}="b_4",
(11,-7)*{}="b_5",
(-3,-17)*{\underbrace{\ \ \ \ \ \ \ \ \ \ \ \ }_A},
\ar @{->} "a";"01" <0pt>
\ar @{->} "a";"02" <0pt>
\ar @{->} "a";"03" <0pt>
\ar @{->} "a";"04" <0pt>
\ar @{<-} "a";"b_2" <0pt>
\ar @{<-} "a";"b_3" <0pt>
\ar @{<-} "a";"b_1" <0pt>
\ar @{<-} "a";"b_4" <0pt>
\ar @{<-} "a";"b_5" <0pt>
\ar @{->} "c1";"b_3" <0pt>
\ar @{->} "c2";"b_3" <0pt>
\ar @{->} "c3";"b_3" <0pt>
\endxy}\Ea
+
\sum_{A\subsetneq [m]\atop \# A\geq 2} 
\Ba{c}\resizebox{22mm}{!}{\xy
(1.6,7)*{...},
(-13,7)*{},
(-8,7)*{},
(-3,7)*{},
(7,7)*{},
(13,7)*{},
 (0,0)*{\bu}="a",
 (0,-6)*{...},
(-5,-7)*{}="01",
(-3,-7)*{}="02",
(3,-7)*{}="03",
(5,-7)*{}="04",
(-12,7)*{}="b_1",
(-9,7)*{}="b_2",
(-6,7)*{...},
(-3,7)*{\bu}="b_3",
(-9,14)*{}="c1",
(-5,14)*{}="c2",
(-1,14)*{...},
(3,14)*{}="c3",
(7,7)*{}="b_4",
(11,7)*{}="b_5",
(-3,17)*{\overbrace{\ \ \ \ \ \ \ \ \ \ \ \ }^A},
\ar @{<-} "a";"01" <0pt>
\ar @{<-} "a";"02" <0pt>
\ar @{<-} "a";"03" <0pt>
\ar @{<-} "a";"04" <0pt>
\ar @{->} "a";"b_2" <0pt>
\ar @{->} "a";"b_3" <0pt>
\ar @{->} "a";"b_1" <0pt>
\ar @{->} "a";"b_4" <0pt>
\ar @{->} "a";"b_5" <0pt>
\ar @{<-} "c1";"b_3" <0pt>
\ar @{<-} "c2";"b_3" <0pt>
\ar @{<-} "c3";"b_3" <0pt>
\endxy}\Ea
.
\Eeq
The properad $\cI\cB_\infty^{path}$ is the associated graded of the dg properad 
$\cI\cB_\infty$ of strongly homotopy infinitesimal bialgebras with respect to the filtration by the total number of directed non-self-intersecting paths connecting in-hairs to out-hairs. The dg properad $\cI\cB_\infty$ was introduced and studied in \cite{A} where it was shown that the cohomology properad $\cI\cB^{path}:=H^\bu(\cI\cB_\infty^{path},\delta_0)$ is generated by the following pair of trivalent ribbon corollas
\Beq\label{3: generators of IB}
 \Ba{c}\resizebox{8mm}{!}{  \xy
(0,6)*{}="1";
    (0,0.2)*{\bu}="L";
  (-4,-5)*{}="C";
   (+4,-5)*{}="D";
\ar @{->} "D";"L" <0pt>
\ar @{->} "C";"L" <0pt>
\ar @{<-} "1";"L" <0pt>
 \endxy}
 \Ea,
 \ \ \
 \Ba{c}\resizebox{8mm}{!}{  \xy
(0,-6)*{}="1";
    (0,-0.2)*{\bu}="L";
  (-4,5)*{}="C";
   (4,5)*{}="D";
\ar @{<-} "D";"L" <0pt>
\ar @{<-} "C";"L" <0pt>
\ar @{->} "1";"L" <0pt>
 \endxy}
 \Ea, 
\Eeq
modulo the following relations
$$
  \Ba{c}\resizebox{6.5mm}{!}{  \xy
(-4,10)*{}="1";
 (4,10)*{}="2";
    (0,3.5)*{\bbu}="A";
 (0,-3.5)*{\bbu}="B";
 (-4,-10)*{}="b1";
 (4,-10)*{}="b2";
\ar @{->} "A";"1" <0pt>
\ar @{->} "A";"2" <0pt>
\ar @{->} "B";"A" <0pt>
\ar @{<-} "B";"b1" <0pt>
\ar @{<-} "B";"b2" <0pt>
 \endxy}
 \Ea
=0,
\ \ \ 
\Ba{c}\resizebox{13.5mm}{!}{  \xy
(-9,10)*{^0}="1";
    (-9,+3)*{\bbu}="L";
 (-14,-3.5)*{\bbu}="B";
 (-20,-12)*+{_1}="b1";
 (-8,-12)*+{_2}="b2";
  (-3,-4)*{_3}="C";
\ar @{->} "C";"L" <0pt>
\ar @{->} "B";"L" <0pt>
\ar @{<-} "B";"b1" <0pt>
\ar @{<-} "B";"b2" <0pt>
\ar @{<-} "1";"L" <0pt>
 \endxy}
 \Ea
 +
  \Ba{c}\resizebox{13.5mm}{!}{  \xy
(9,10)*{^0}="1";
    (9,+3)*{\bbu}="L";
 (14,-3.5)*{\bbu}="B";
 (20,-12)*+{_3}="b1";
 (8,-12)*+{2}="b2";
  (3,-4)*{1}="C";
\ar @{->} "C";"L" <0pt>
\ar @{->} "B";"L" <0pt>
\ar @{<-} "B";"b1" <0pt>
\ar @{<-} "B";"b2" <0pt>
\ar @{<-} "1";"L" <0pt>
 \endxy}
 \Ea=0, 
 \Ba{c}\resizebox{13.5mm}{!}{  \xy
(-9,-10)*{^0}="1";
    (-9,-3)*{\bbu}="L";
 (-14,3.5)*{\bbu}="B";
 (-20,12)*+{_1}="b1";
 (-8,12)*+{_2}="b2";
  (-3,4)*{_3}="C";
\ar @{->} "C";"L" <0pt>
\ar @{->} "B";"L" <0pt>
\ar @{<-} "B";"b1" <0pt>
\ar @{<-} "B";"b2" <0pt>
\ar @{<-} "1";"L" <0pt>
 \endxy}
 \Ea
 +
  \Ba{c}\resizebox{13.5mm}{!}{  \xy
(9,-10)*{^0}="1";
    (9,-3)*{\bbu}="L";
 (14,3.5)*{\bbu}="B";
 (20,12)*+{_3}="b1";
 (8,12)*+{2}="b2";
  (3,4)*{1}="C";
\ar @{<-} "C";"L" <0pt>
\ar @{<-} "B";"L" <0pt>
\ar @{->} "B";"b1" <0pt>
\ar @{->} "B";"b2" <0pt>
\ar @{->} "1";"L" <0pt>
 \endxy}
 \Ea=0.
$$

\subsubsection{\sf Remark on boundedness of filtrations}
It is important to notice that the generators of both complexes are hairy ribbon graphs for which the genus $g$ and of the number $b$ of boundaries are well-defined; moreover these integer parameters  $g$ and $b$ are preserved by the differential in $\PCYf_d(m,n)$ and by relations (\ref{2: relations I in DeltaPCYf})-(\ref{2: double Lie relations}). Hence both complexes decomposed into the direct sums of subcomplexes
$$
\PCYf_d(m,n)=\bigoplus_{g,b}\PCYf_d(m,n)^{g,b}, \ \ \ \ 
\Delta\PCYf_d(m,n)=\bigoplus_{g,b}\PCYf_d(m,n)^{g,b}
$$
parameterized by the genus and the number of boundaries of the generators.
The Main Theorem is proven once we show an isomorphism of the cohomology groups 
$$
H^\bu(\PCYf_d(m,n)^{g,b})=H^\bu(\Delta\PCYf_d(m,n)^{g,b}).
$$
It is easy to see that for any fixed values of the integer parameters $m,n,g,b$ the subspaces of both complexes spanned by hairy ribbon graph
with any fixed cohomological degree $k$ are {\it finite-dimensional}. Hence the spectral sequences of the filtrations of both complexes which we consider below must be convergent.

\subsubsection{\sf Filtration by the number of sources}
Consider filtrations of both sides of the epimorphism
$$
\pi: \PCYf_d\lon \Delta \PCYf_d
$$
by the number of sources, that is, vertices with no ingoing edge or hair
as in the picture (\ref{1: sources}). The Main Theorem is proven once we show that the induced morphism
\Beq\label{2: gr(p_1)}
gr(\pi): gr(\PCYf_d)\lon gr(\Delta \PCYf_d)
\Eeq
of the associated graded complexes is a quasi-isomorphism. Note that the induced differential  $\delta=\sum_v\delta_v$ acts trivially on both bivalent and trivalent vertices $v$ in these associated graded complexes.

\sip

Let  $\PCYf_d^{marked}$ and $\Delta \PCYf_d^{marked}$ be the versions of
$gr(\PCYf_d)$ and, respectively, $gr(\Delta \PCYf_d)$ in which sources are distinguished, say labelled by some set of integers. There is an associated epimorphism of complexes
\Beq\label{2: map pi-marked}
\pi^{marked}: \PCYf_d^{marked}\lon \Delta \PCYf_d^{marked}.
\Eeq
By Maschke's  theorem, the Main Theorem is proven once we show that the epimorphism $\pi^{marked}$ is a quasi-isomorphism.

\subsubsection{\sf An intermediate complex} A
vertex $v$ of valency $\geq 4$
can have  two or more neighboring (with respect to the given cyclic order) edges having the same direction; we call such edges {\it parallel at $v$}. For example, the vertex
$$
v=  \Ba{c}\resizebox{13mm}{!}
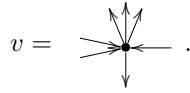
{  \xy
(0,7)*{}="U";
(-2.5,6)*{}="UL";
(2.5,6)*{}="UR";
(0,-6)*{}="D";
    (0,0)*{\bu}="C";
(-7,1.5)*{}="L1";
(-7,-1.5)*{}="L2";    
  (7,0)*{}="R";
  
\ar @{->} "C";"U" <0pt>
\ar @{->} "C";"UL" <0pt>
\ar @{->} "C";"UR" <0pt>
\ar @{->} "C";"D" <0pt>
\ar @{<-} "C";"L1" <0pt>
\ar @{<-} "C";"L2" <0pt>
\ar @{<-} "C";"R" <0pt>
 \endxy}
 \Ea.
$$
 has three parallel out-edges and two parallel in-edges. Parallel edges come in maximal cyclically ordered blocks which are called {\it bunches}\, at $v$. The above vertex has 4 different bunches.

 Let $K$ be a linear subspace 
of the complex $\PCYf_d^{marked}$ spanned by hairy ribbon graphs having at least one source of valency $\geq 3$ or at least one vertex $v$ of valency $\geq 4$ which has at least one pair of parallel edges at $v$, and let $\langle K, \delta K\rangle$ be its differential closure in $\PCYf_d^{marked}$. Define an intermediate  complex  $\widehat{\Delta}\PCYf_d^{marked}$  by the following short exact sequence
$$
0 \lon \langle K, \delta K\rangle \lon \PCYf_d^{marked} \stackrel{p_1}{\lon} \widehat{\Delta} \PCYf_d^{marked} \lon 0.
$$
The quotient complex  $\widehat{\Delta} \PCYf_d^{marked}$ is generated by equivalence classes of hairy ribbon graphs which can have 2- and 3-valent vertices of the form
\Beq\label{2: 2 and 3 valent vertices}
\Ba{c}\resizebox{8mm}{!}{ \xy
 (0,0)*{}="a",
(4,-3)*{\bullet}="b",
(8,0)*{}="c",
\ar @{<-} "a";"b" <0pt>
\ar @{->} "b";"c" <0pt>
\endxy}\Ea, 
\ \ \ 
 \Ba{c}\resizebox{8mm}{!}{  \xy
(0,6)*{}="1";
    (0,0.2)*{\bu}="L";
  (-4,-5)*{}="C";
   (+4,-5)*{}="D";
\ar @{->} "D";"L" <0pt>
\ar @{->} "C";"L" <0pt>
\ar @{<-} "1";"L" <0pt>
 \endxy}
 \Ea,
 \ \ \
 \Ba{c}\resizebox{8mm}{!}{  \xy
(0,-6)*{}="1";
    (0,-0.2)*{\bu}="L";
  (-4,5)*{}="C";
   (4,5)*{}="D";
\ar @{<-} "D";"L" <0pt>
\ar @{<-} "C";"L" <0pt>
\ar @{->} "1";"L" <0pt>
 \endxy}
 \Ea, 
\Eeq
and also vertices of valencies $2k$, $k\geq 2$, with no parallel edges/hairs attached
as e.g.\ the following ones
\Beq\label{2: DLie corollas}
  \Ba{c}\resizebox{12mm}{!}{  \xy
(0,-6)*{}="U";
(0,+6)*{}="D";
    (0,0)*{\bu}="C";
  (6,0)*{}="R";
   (-6,0)*{}="L";
\ar @{->} "C";"U" <0pt>
\ar @{->} "C";"D" <0pt>
\ar @{<-} "C";"L" <0pt>
\ar @{<-} "C";"R" <0pt>
 \endxy}
 \Ea
 ,\
   \Ba{c}\resizebox{12mm}{!}{  \xy
   (-4,6)*{}="UL";
(4,6)*{}="UR";
(-4,-6)*{}="DL";
(4,-6)*{}="DR";
    (0,0)*{\bu}="C";
  (6,0)*{}="R";
   (-6,0)*{}="L";
\ar @{->} "C";"UL" <0pt>
\ar @{<-} "C";"UR" <0pt>
\ar @{->} "C";"DL" <0pt>
\ar @{<-} "C";"DR" <0pt>
\ar @{<-} "C";"L" <0pt>
\ar @{->} "C";"R" <0pt>
 \endxy}
 \Ea.
\Eeq
The induced differential $\delta$ in $\widehat{\Delta}  \PCYf_d^{marked}$ acts  non-trivially only  on vertices with valency $\geq 4$.

\sip

It is worth noting the subspace $L$ of  $\widehat{\Delta}  \PCYf_d^{marked}$ spanned by ribbon graphs having at least one vertex of valency $\leq 3$ is a subcomplex such that the quotient complex
$\widehat{\Delta}  \PCYf_d^{marked}/L$ is precisely the minimal resolution $\caD\Lie_{d,\infty}$ of the Koszul properad $\caD\Lie_d$ of double Lie algebras which is generated by the 4-valent corolla  
$  \Ba{c}\resizebox{10mm}{!}{  \xy
(0,-5)*{}="U";
(0,+5)*{}="D";
    (0,0)*{\bu}="C";
  (5,0)*{}="R";
   (-5,0)*{}="L";
\ar @{->} "C";"U" <0pt>
\ar @{->} "C";"D" <0pt>
\ar @{<-} "C";"L" <0pt>
\ar @{<-} "C";"R" <0pt>
 \endxy}
 \Ea$ modulo the relation (\ref{2: double Lie relations}). Koszulness of $\caD\Lie_d$  was proven by J.\ Leray in \cite{L}.
Thanks to this result we have the equality,
$$
H^\bu(\widehat{\Delta}  \PCYf_d^{marked}/L)=\caD\Lie_{d},
$$
which plays an important role in our arguments below.

\sip

\subsubsection{\sf Claim I}\label{2: Claim 1} {\it  The epimorphism $p_1: \PCYf_d^{marked} {\lon} \widehat{\Delta} \PCYf_d^{marked}$ is a quasi-isomorphism}.

\mip

\noindent Let us prove this claim by considering a filtration of both sides of the arrow $p_1$ by the total number of directed paths  connecting the sources and the incoming hairs of ribbon graphs  to their outgoing hairs. There is an induced morphism of associated graded complexes
$$
gr(p_1):  gr(\PCYf_d^{marked}) {\lon} gr(\widehat{\Delta} \PCYf_d^{marked}).
$$
The Claim is proven once we show that $gr(p_1)$ is a quasi-isomorphism.

\sip

The complex $gr(\widehat{\Delta} \PCYf_d^{marked})$ has a trivial differential. The relations (\ref{2: relations I in DeltaPCYf}) and (\ref{2: 3+4 relations}) imposed on the generators of $\widehat{\Delta} \PCYf_d^{marked}$ simplify in the associated graded complex to the following ones
\Beq\label{2: simplified relations}
  \Ba{c}\resizebox{5.0mm}{!}{  \xy
(-4,10)*{}="1";
 (4,10)*{}="2";
    (0,3.5)*{\bbu}="A";
 (0,-3.5)*{\bbu}="B";
 (-4,-10)*{}="b1";
 (4,-10)*{}="b2";
\ar @{->} "A";"1" <0pt>
\ar @{->} "A";"2" <0pt>
\ar @{->} "B";"A" <0pt>
\ar @{<-} "B";"b1" <0pt>
\ar @{<-} "B";"b2" <0pt>
 \endxy}
 \Ea
\hspace{-2mm} =0,
\Ba{c}\resizebox{10.5mm}{!}{  \xy
(-9,10)*{^0}="1";
    (-9,+3)*{\bbu}="L";
 (-14,-3.5)*{\bbu}="B";
 (-20,-12)*+{_1}="b1";
 (-8,-12)*+{_2}="b2";
  (-3,-4)*{_3}="C";
\ar @{->} "C";"L" <0pt>
\ar @{->} "B";"L" <0pt>
\ar @{<-} "B";"b1" <0pt>
\ar @{<-} "B";"b2" <0pt>
\ar @{<-} "1";"L" <0pt>
 \endxy}
 \Ea
 +
  \Ba{c}\resizebox{10.5mm}{!}{  \xy
(9,10)*{^0}="1";
    (9,+3)*{\bbu}="L";
 (14,-3.5)*{\bbu}="B";
 (20,-12)*+{_3}="b1";
 (8,-12)*+{2}="b2";
  (3,-4)*{1}="C";
\ar @{->} "C";"L" <0pt>
\ar @{->} "B";"L" <0pt>
\ar @{<-} "B";"b1" <0pt>
\ar @{<-} "B";"b2" <0pt>
\ar @{<-} "1";"L" <0pt>
 \endxy}
 \Ea
 \hspace{-2mm}=0, 
 \Ba{c}\resizebox{10.5mm}{!}{  \xy
(-9,-10)*{^0}="1";
    (-9,-3)*{\bbu}="L";
 (-14,3.5)*{\bbu}="B";
 (-20,12)*+{_1}="b1";
 (-8,12)*+{_2}="b2";
  (-3,4)*{_3}="C";
\ar @{->} "C";"L" <0pt>
\ar @{->} "B";"L" <0pt>
\ar @{<-} "B";"b1" <0pt>
\ar @{<-} "B";"b2" <0pt>
\ar @{<-} "1";"L" <0pt>
 \endxy}
 \Ea
 +
  \Ba{c}\resizebox{10.5mm}{!}{  \xy
(9,-10)*{^0}="1";
    (9,-3)*{\bbu}="L";
 (14,3.5)*{\bbu}="B";
 (20,12)*+{_3}="b1";
 (8,12)*+{2}="b2";
  (3,4)*{1}="C";
\ar @{<-} "C";"L" <0pt>
\ar @{<-} "B";"L" <0pt>
\ar @{->} "B";"b1" <0pt>
\ar @{->} "B";"b2" <0pt>
\ar @{->} "1";"L" <0pt>
 \endxy}
 \Ea
 \hspace{-2mm}=0,
   \Ba{c}\resizebox{16mm}{!}{  \xy
(-3,13)*+{_1}="UL";
(3,13)*+{_2}="UR";
(0,7)*{\bu}="U";
(0,-7)*+{_{2k}}="D";
(0,-10)*{\underbrace{\ \ \ \ \ \ \ \ \ \ \ \ \ \ \ \ \ }_{2k\ \text{valent vertex}}};
    (0,0)*{\bu}="C";
  (7,0)*+{_3}="R";
   (-7,0)*+{_0}="L";
\ar @{->} "U";"UL" <0pt>   
\ar @{->} "U";"UR" <0pt>  
\ar @{->} "C";"U" <0pt>
\ar @{->} "C";"D" <0pt>
\ar @{<-} "C";"L" <0pt>
\ar @{<-} "C";"R" <0pt>
 \endxy}
 \Ea
\hspace{-2mm}=0,
   \Ba{c}\resizebox{16mm}{!}{  \xy
(-3,13)*+{_1}="UL";
(3,13)*+{_2}="UR";
(0,7)*{\bu}="U";
(0,-7)*+{_{2k}}="D";
(0,-10)*{\underbrace{\ \ \ \ \ \ \ \ \ \ \ \ \ \ \ \ \ }_{2k\ \text{valent vertex}}};
    (0,0)*{\bu}="C";
  (7,0)*+{_3}="R";
   (-7,0)*+{_0}="L";
\ar @{<-} "U";"UL" <0pt>   
\ar @{<-} "U";"UR" <0pt>  
\ar @{<-} "C";"U" <0pt>
\ar @{<-} "C";"D" <0pt>
\ar @{->} "C";"L" <0pt>
\ar @{->} "C";"R" <0pt>
 \endxy}
 \Ea\hspace{-2mm}=0.
\Eeq
This relations can be easily solved so that  $\widehat{\Delta} \PCYf_d^{marked}$ gets a nice presentation in terms of decorated {\it reduced}\, ribbon graphs which are described below.

\sip

The differential $\delta_0$ on the complex $gr(\PCYf_d^{marked})$ acts non-trivially only on sources (\ref{2: sources}) of valency $\geq 3$ and on those vertices of valency $\geq 4$ which have at least one pair of parallel edges. The edges attached to any vertex $v$ of a generator $\Ga\in gr(\PCYf_d^{marked})$ can be grouped in bunches consisting of {\it all}\, neighboring edges parallel to each other; the cardinality of such a bunch is $\geq 1$, and the total number $b(v)$ of such bunches at $v$  is equal to $1$ if $v$ is a source, or it is equal to an even number otherwise.
 For example, the vertices in the following set
$$
v_1=  \Ba{c}\resizebox{13mm}{!}{  \xy
(0,-5)*{}="U";
(0,+5)*{}="D";
    (0,0)*{\bu}="C";
  (5,0)*{}="R";
   (-5,0)*{}="L";
\ar @{->} "C";"U" <0pt>
\ar @{->} "C";"D" <0pt>
\ar @{<-} "C";"L" <0pt>
\ar @{<-} "C";"R" <0pt>
 \endxy}
 \Ea, \ \ 
v_2=  \Ba{c}\resizebox{13mm}{!}{  \xy
(0,5)*{}="U";
(-2,4)*{}="UL";
(2,4)*{}="UR";
(0,-5)*{}="D";
    (0,0)*{\bu}="C";
(-5,1)*{}="L1";
(-5,-1)*{}="L2";    
  (5,0)*{}="R";
  
\ar @{->} "C";"U" <0pt>
\ar @{->} "C";"UL" <0pt>
\ar @{->} "C";"UR" <0pt>
\ar @{->} "C";"D" <0pt>
\ar @{<-} "C";"L1" <0pt>
\ar @{<-} "C";"L2" <0pt>
\ar @{<-} "C";"R" <0pt>
 \endxy}
 \Ea
 ,
 \ \ 
 v_3=
 \Ba{c}\resizebox{11mm}{!}{ \xy
(0,8)*{\overbrace{\  \ \ \ \ \ \ \ \ \ }^{m\geq 2}},
(0,-8)*{\underbrace{\  \ \ \ \ \ \ \ \ \ }_{n\geq 2}},
(0,4.5)*+{...},
(0,-4.5)*+{...},
(0,0)*{\bu}="o",
(-5,5)*{}="1",
(-3,5)*{}="2",
(3,5)*{}="3",
(5,5)*{}="4",
(-3,-5)*{}="5",
(3,-5)*{}="6",
(5,-5)*{}="7",
(-5,-5)*{}="8",
\ar @{->} "o";"1" <0pt>
\ar @{->} "o";"2" <0pt>
\ar @{->} "o";"3" <0pt>
\ar @{->} "o";"4" <0pt>
\ar @{<-} "o";"5" <0pt>
\ar @{<-} "o";"6" <0pt>
\ar @{<-} "o";"7" <0pt>
\ar @{<-} "o";"8" <0pt>
\endxy}\Ea
$$
have the number of bunches given by
 $$
 b(v_1)=b(v_2)=4,\ \ \ b(v_3)=2.
 $$
 Let us call vertices $v$ of $\Ga$ {\it operadic}\, if
$b(v)=2$ and precisely one bunch at $v$ has cardinality 1. Thus operadic vertices have one of the following two forms
$$
  \underbrace{\Ba{c}\resizebox{20mm}{!}
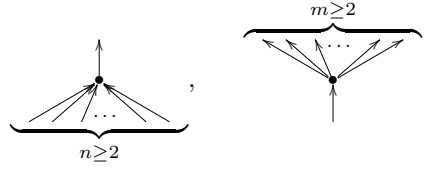
{\xy
(1,-6)*{\ldots},
(-13,-7)*{},
(-8,-7)*{},
(-3,-7)*{},
(7,-7)*{},
(13,-7)*{},
 (0,0)*{\bu}="a",
(0,7)*{}="0",
(-12,-7)*{}="b_1",
(-8,-7)*{}="b_2",
(-3,-7)*{}="b_3",
(8,-7)*{}="b_4",
(12,-7)*{}="b_5",
\ar @{->} "a";"0" <0pt>
\ar @{<-} "a";"b_2" <0pt>
\ar @{<-} "a";"b_3" <0pt>
\ar @{<-} "a";"b_1" <0pt>
\ar @{<-} "a";"b_4" <0pt>
\ar @{<-} "a";"b_5" <0pt>
\endxy}\Ea}_{n\geq 2}, \ \ \ \ \ 
  \overbrace{\Ba{c}\resizebox{20mm}{!}{\xy
(1,6)*{\ldots},
(-13,7)*{},
(-8,7)*{},
(-3,7)*{},
(7,7)*{},
(13,7)*{},
 (0,0)*{\bu}="a",
(0,-7)*{}="0",
(-12,7)*{}="b_1",
(-8,7)*{}="b_2",
(-3,7)*{}="b_3",
(8,7)*{}="b_4",
(12,7)*{}="b_5",
\ar @{<-} "a";"0" <0pt>
\ar @{->} "a";"b_2" <0pt>
\ar @{->} "a";"b_3" <0pt>
\ar @{->} "a";"b_1" <0pt>
\ar @{->} "a";"b_4" <0pt>
\ar @{->} "a";"b_5" <0pt>
\endxy}\Ea}^{m\geq 2}.
$$
Operadic vertices have a distinguished edge which belongs to the bunch
of cardinality one; let us call such an edge {\it operadic}\, as well.
As the induced differential $\delta_0$  must preserve the path filtration, it can create only operadic vertices as shown, for example, in (\ref{2: d in Ass_infty}), (\ref{2: d in Ass_wedge_infty}),(\ref{2: d in IB_infty_path}) and in the following example, 
$$
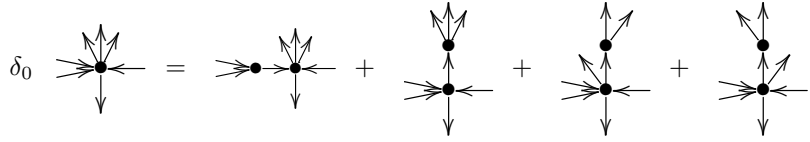
\delta_0 \Ba{c}\resizebox{13mm}{!}{  \xy
(0,5)*{}="U";
(-2,4)*{}="UL";
(2,4)*{}="UR";
(0,-5)*{}="D";
    (0,0)*{\bu}="C";
(-5,1)*{}="L1";
(-5,-1)*{}="L2";    
  (5,0)*{}="R";
  
\ar @{->} "C";"U" <0pt>
\ar @{->} "C";"UL" <0pt>
\ar @{->} "C";"UR" <0pt>
\ar @{->} "C";"D" <0pt>
\ar @{<-} "C";"L1" <0pt>
\ar @{<-} "C";"L2" <0pt>
\ar @{<-} "C";"R" <0pt>
 \endxy}
 \Ea
 =
  \Ba{c}\resizebox{17mm}{!}{  \xy
(0,5)*{}="U";
(-2,4)*{}="UL";
(2,4)*{}="UR";
(0,-5)*{}="D";
    (0,0)*{\bu}="C";
     (-5,0)*{\bu}="CL";
(-10,1)*{}="L1";
(-10,-1)*{}="L2";    
  (5,0)*{}="R";
  \ar @{->} "CL";"C" <0pt>
\ar @{->} "C";"U" <0pt>
\ar @{->} "C";"UL" <0pt>
\ar @{->} "C";"UR" <0pt>
\ar @{->} "C";"D" <0pt>
\ar @{<-} "CL";"L1" <0pt>
\ar @{<-} "CL";"L2" <0pt>
\ar @{<-} "C";"R" <0pt>
 \endxy}
 \Ea
 +
  \Ba{c}\resizebox{13mm}{!}{  \xy
(0,10)*{}="U";
(-2,9)*{}="UL";
(2,9)*{}="UR";
(0,-5)*{}="D";
    (0,0)*{\bu}="C";
      (0,5)*{\bu}="CU";
(-5,1)*{}="L1";
(-5,-1)*{}="L2";    
  (5,0)*{}="R";
  \ar @{->} "C";"CU" <0pt>
\ar @{->} "CU";"U" <0pt>
\ar @{->} "CU";"UL" <0pt>
\ar @{->} "CU";"UR" <0pt>
\ar @{->} "C";"D" <0pt>
\ar @{<-} "C";"L1" <0pt>
\ar @{<-} "C";"L2" <0pt>
\ar @{<-} "C";"R" <0pt>
 \endxy}
 \Ea
 +
   \Ba{c}\resizebox{13mm}{!}{  \xy
(0,10)*{}="U";
(-3,4)*{}="UL";
(3,9)*{}="UR";
(0,-5)*{}="D";
    (0,0)*{\bu}="C";
      (0,5)*{\bu}="CU";
(-5,1)*{}="L1";
(-5,-1)*{}="L2";    
  (5,0)*{}="R";
  \ar @{->} "C";"CU" <0pt>
\ar @{->} "CU";"U" <0pt>
\ar @{->} "C";"UL" <0pt>
\ar @{->} "CU";"UR" <0pt>
\ar @{->} "C";"D" <0pt>
\ar @{<-} "C";"L1" <0pt>
\ar @{<-} "C";"L2" <0pt>
\ar @{<-} "C";"R" <0pt>
 \endxy}
 \Ea
  +
   \Ba{c}\resizebox{13mm}{!}{  \xy
(0,10)*{}="U";
(-3,9)*{}="UL";
(3,4)*{}="UR";
(0,-5)*{}="D";
    (0,0)*{\bu}="C";
      (0,5)*{\bu}="CU";
(-5,1)*{}="L1";
(-5,-1)*{}="L2";    
  (5,0)*{}="R";
  \ar @{->} "C";"CU" <0pt>
\ar @{->} "CU";"U" <0pt>
\ar @{->} "CU";"UL" <0pt>
\ar @{->} "C";"UR" <0pt>
\ar @{->} "C";"D" <0pt>
\ar @{<-} "C";"L1" <0pt>
\ar @{<-} "C";"L2" <0pt>
\ar @{<-} "C";"R" <0pt>
 \endxy}
 \Ea
$$
Note that each bunch of a vertex $v$ with $b(v)\geq 3$ is controlled essentially by the complex $\cA ss_\infty^{\wedge,+}$ or the complex $\cA ss_\infty^{\vee,+}$ depending on the direction of edges in that bunch.

\sip

 Next we decompose our complex $(gr(\PCYf_d^{marked}), \delta_0)$ into the tensor product of complexes parameterized by so-called {\it reduced}\, ribbon graphs with hairs defined as follows:
 given any generator $\Ga\in gr(\PCYf_d^{marked})$, let $\widetilde{\Ga}$ be the so-called {\it reduced} ribbon graph obtained from $\Ga$ by contracting all its operadic internal edges; note that this procedure  does {\it not}\, create {\it closed}\, paths of directed internal edges in $\widetilde{\Ga}$.  Denote the set of such ribbon graphs with hairs without operadic edges by $\widetilde{S}_{red}$.

\sip

For any $G\in \widetilde{S}_{red}$ denote by  $gr_G(\PCYf_d^{marked})$
the linear subspace of $gr_G(\PCYf_d^{marked})$ spanned by all hairy ribbon graphs
$\Ga$ satisfying the condition $\widetilde{\Ga}=G$. As $\widetilde{(\delta_0\Ga)}=\widetilde{\Ga}$, this linear subspace is a subcomplex of $(gr(\PCYf_d^{marked}), \delta_0)$; moreover we have a direct sum decomposition
$$
(gr(\PCYf_d^{marked}), \delta_0)=  \bigoplus_{G\in \widetilde{S}_{red}} (gr_G(\PCYf_d^{marked}), \delta_0).
$$
Similarly one  has a direct sum decomposition of the trivial complex,
$$
gr(\widehat{\Delta} \PCYf_d^{marked})= \bigoplus_{G\in \widetilde{S}_{red}} gr_G(\widehat{\Delta} \PCYf_d^{marked})
$$ 
The Claim {\ref{2: Claim 1}} is finally proven once we show that for any $G\in \widetilde{S}_{red}$ one has an isomorphism of graded vector spaces 
\Beq\label{2: gr G isomorphism}
 H^\bu(gr_G(\PCYf_d^{marked}), \delta_0) \simeq gr_G(\widehat{\Delta} \PCYf_d^{marked})
\Eeq
as this implies that the epimorphism $gr(p_1)$ in (\ref{2: gr(p_1)}) is a quasi-isomorphism.

\sip

Let us prove the isomorphism (\ref{2: gr G isomorphism}) for arbitrary reduced graph $G\in  \widetilde{S}_{red}$. Denote by
\Bi
\item[(1)] $V_s(G)$ the set of sources of $G$; the valency of $v\in V_s(G)$ is denoted by $|v|_{out}$.

\item[(2)] $V_2(G)$ the set of vertices with $b(v)=2$; such a vertex has $|v|_{in}\geq 2$ incoming edges and  $|v|_{out}\geq 2$  outgoing edges;
\item[(3)] $V_{>2}(G)$ the set of vertices with $b(v)>2$; the set of bunches at $v$ with out-going (resp., incoming)  edges is denoted by $B_{out}(v)$ (resp.,  $B_{in}(v)$),
    and the cardinality of a bunch $b$ in   $B_{out}(v)$ (resp., in  $B_{in}(v)$) is denoted by $|b|$.
\Ei
Using this notation we can decompose the complex $(gr_G(\PCYf_d^{marked}), \delta_0)$ into the unordered tensor product of complexes parameterized by vertices $v$ of the fixed reduced ribbon graph $G$,
\Beqrn
gr_G(\PCYf_d^{marked})\hspace{-2mm}& \simeq & \hspace{-2mm}\left( \bigotimes_{v\in V_s(G)}  \cA ss_\infty^{\vee,cyc}(|v|_{out})\right) \ot \left( \bigotimes_{v\in V_2(G)} gr \cI\cB_\infty (|v|_{out}, |v|_{in})\right) \\
&& \ot \left(\bigotimes_{v\in V_{>2}(G)}\left( \bigotimes_{b\in B_{in}(v)}  \cA ss_\infty^{\wedge,+}(|b|)  \bigotimes_{b\in B_{out}(v)}  \cA ss_\infty^{\vee,+}(|b|)        \right)\right)
\Eeqrn
Using now the well-known results (see \S {\ref{2: subsec on aux complexes}}) on cohomology groups of the above tensor factors, one obtains the required isomorphism (\ref{2: gr G isomorphism}). The Claim {\ref{2: Claim 1}} is proven.

\newcommand{\sRGC}{{\Delta} \PCYf^{marked}}

\subsubsection{\sf From the intermediate complex to $\Delta \PCYf_d^{marked}$}
Let $\langle L, \delta L\rangle$  be the differential closure in the intermediate complex $\widehat{\Delta} \PCYf_d^{marked}$ of the linear subspace $L$ spanned by (equivalence classes of) hairy ribbon graphs  having at least one 
vertex $v$ of valency $>4$. There is  short exact sequence of complexes
$$
0 \lon \langle L, \delta L\rangle \lon \widehat{\Delta} \PCYf_d^{marked}\stackrel{p_2}{\lon} 
\sRGC_d \lon 0.
$$
The epimorphism of complexes (\ref{2: map pi-marked}) factors through the composition
$$
\pi^{marked}: \PCYf_d^{marked}\stackrel{p_1}{\lon} \wh{\Delta} \PCYf_d^{marked}
\stackrel{p_2}{\lon} \Delta\PCYf_d^{marked}.
$$
We have shown above that $p_1$ is a quasi-isomorphism. Thus to show that
 $\pi^{marked}$ is a quasi-isomorphism (and hence complete the proof of Main Theorem {\ref{1: Theorem on PCYf}}) it remains to prove the following

\sip

\noindent{\sf Claim II}: {\it The epimorphism $\widehat{\Delta}\PCYf_d^{marked}\stackrel{p_2}{\lon} 
\sRGC_d$ is a quasi-isomorphism.} 

\begin{proof}
 Consider a  filtration 
 $$
 F_0\subset \ldots \subset F_p \subset F_{p+1}\subset\ldots
 $$
  of both sides of the epimorphism $p_2$  by the number of 3-valent vertices, and let
$$
Gr(p_2): Gr(\widehat{\Delta}\PCYf_d^{marked}) \lon Gr (\sRGC_d)
$$
be the induced morphism  of the associated graded. If we show that $Gr(p_2)$
is a quasi-isomorphism, the Claim II is proven.
The complex $Gr (\sRGC_2)$ has trivial differential so that it is enough to show that 
\Beq\label{2: final iso for wh(Delta)}
H^\bu(Gr(\widehat{\Delta}\PCYf_d^{marked})) \simeq Gr (\sRGC_d).
\Eeq

The induced differential $d$ on the complex $Gr(\widehat{\Delta}\PCYf_d^{marked})$ acts
non-trivially only on vertices of  valencies $>4$, splitting them precisely as the differential in the dg free properad $\caD\Lie_\infty$  strongly homotopy double Lie algebras. The latter properad is freely generated by at least 4-valent  
 ribbon corollas with no parallel hairs as in (\ref{2: DLie corollas}), and the differential $d$ on $\caD\Lie_\infty$ acts on generators by substituting into the vertex the graph $\xy
 (0,0)*{\bu}="a",
(5,0)*{\bu}="b",
\ar @{->} "a";"b" <0pt>
\endxy$ and then taking a sum over all possible ways to reattach edges and hairs (attached earlier to $v$) among the two newly created vertices while (i) respecting the natural cyclic structure on that edges and hairs, and (ii) setting to zero terms which have a vertex of valency $\leq 3$ or a vertex with at least of pair of parallel edges/hairs, for example
\Beqrn
d\left(  \Ba{c}\resizebox{13mm}{!}{  \xy
(0,-5)*{}="U";
(0,+5)*{}="D";
    (0,0)*{\bu}="C";
  (5,0)*{}="R";
   (-5,0)*{}="L";
\ar @{->} "C";"U" <0pt>
\ar @{->} "C";"D" <0pt>
\ar @{<-} "C";"L" <0pt>
\ar @{<-} "C";"R" <0pt>
 \endxy}
 \Ea\right)&=&0,\\
d\left( \Ba{c}\resizebox{16mm}{!}{  \xy
   (-5,6)*{_1}="UL";
(5,6)*{_2}="UR";
(-5,-6)*{_5}="DL";
(5,-6)*{_4}="DR";
    (0,0)*{\bu}="C";
  (7,0)*{_3}="R";
   (-7,0)*{_0}="L";
\ar @{->} "C";"UL" <0pt>
\ar @{<-} "C";"UR" <0pt>
\ar @{->} "C";"DL" <0pt>
\ar @{<-} "C";"DR" <0pt>
\ar @{<-} "C";"L" <0pt>
\ar @{->} "C";"R" <0pt>
 \endxy}
 \Ea\right)
& = &
   \Ba{c}\resizebox{28mm}{!}{  \xy
(0,-7)*+{_5}="D";
(0,+7)*+{_1}="U";
 (-7,0)*+{_0}="L";
(0,0)*{\bu}="C";
(7,0)*{\bu}="CR";
  (7,-7)*+{_4}="RD";
(7,+7)*+{_2}="RU";
 (14,0)*+{_3}="RR";
\ar @{->} "C";"U" <0pt>
\ar @{->} "C";"D" <0pt>
\ar @{<-} "C";"L" <0pt>
\ar @{<-} "C";"CR" <0pt>
\ar @{<-} "CR";"RU" <0pt>
\ar @{<-} "CR";"RD" <0pt>
\ar @{->} "CR";"RR" <0pt>
 \endxy}
 \Ea
 +
  \Ba{c}\resizebox{28mm}{!}{  \xy
(0,-7)*+{_3}="D";
(0,+7)*+{_5}="U";
 (-7,0)*+{_4}="L";
(0,0)*{\bu}="C";
(7,0)*{\bu}="CR";
  (7,-7)*+{_2}="RD";
(7,+7)*+{_0}="RU";
 (14,0)*+{_1}="RR";
\ar @{->} "C";"U" <0pt>
\ar @{->} "C";"D" <0pt>
\ar @{<-} "C";"L" <0pt>
\ar @{<-} "C";"CR" <0pt>
\ar @{<-} "CR";"RU" <0pt>
\ar @{<-} "CR";"RD" <0pt>
\ar @{->} "CR";"RR" <0pt>
 \endxy}
 \Ea
  +
  \Ba{c}\resizebox{28mm}{!}{  \xy
(0,-7)*+{_1}="D";
(0,+7)*+{_3}="U";
 (-7,0)*+{_2}="L";
(0,0)*{\bu}="C";
(7,0)*{\bu}="CR";
  (7,-7)*+{_0}="RD";
(7,+7)*+{_4}="RU";
 (14,0)*+{_5}="RR";
\ar @{->} "C";"U" <0pt>
\ar @{->} "C";"D" <0pt>
\ar @{<-} "C";"L" <0pt>
\ar @{<-} "C";"CR" <0pt>
\ar @{<-} "CR";"RU" <0pt>
\ar @{<-} "CR";"RD" <0pt>
\ar @{->} "CR";"RR" <0pt>
 \endxy}
 \Ea.
\Eeqrn  
 Since the numbers of bivalent and trivalent vertices in the complex $Gr(\wh{\Delta}\PCYf_d^{marked})$ are fixed, we can assume without loss of generality that they all are distinguished. Moreover, we can in this case easily solve almost all of the relations (\ref{2: relations I in DeltaPCYf})-(\ref{2: 3+4 relations}) in  $Gr(\wh{\Delta}\PCYf_d^{marked})$ and assume from now on that it is spanned by hairy ribbon graph with {\it no internal edges}\, as in the following list,
 $$
 \Ba{c}\resizebox{5.0mm}{!}
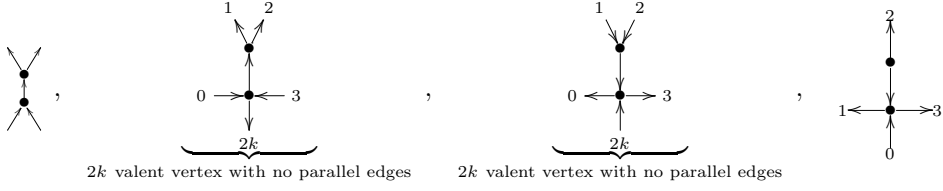
{  \xy
(-4,10)*{}="1";
 (4,10)*{}="2";
    (0,3.5)*{\bbu}="A";
 (0,-3.5)*{\bbu}="B";
 (-4,-10)*{}="b1";
 (4,-10)*{}="b2";
\ar @{->} "A";"1" <0pt>
\ar @{->} "A";"2" <0pt>
\ar @{->} "B";"A" <0pt>
\ar @{<-} "B";"b1" <0pt>
\ar @{<-} "B";"b2" <0pt>
 \endxy}
 \Ea,
    \Ba{c}\resizebox{44mm}{!}{  \xy
(-3,13)*+{_1}="UL";
(3,13)*+{_2}="UR";
(0,7)*{\bu}="U";
(0,-7)*+{_{2k}}="D";
(0,-10)*{\underbrace{\ \ \ \ \ \ \ \ \ \ \ \ \ \ \ \ \ }_{2k\ \text{valent vertex with no parallel edges}}};
    (0,0)*{\bu}="C";
  (7,0)*+{_3}="R";
   (-7,0)*+{_0}="L";
\ar @{->} "U";"UL" <0pt>   
\ar @{->} "U";"UR" <0pt>  
\ar @{->} "C";"U" <0pt>
\ar @{->} "C";"D" <0pt>
\ar @{<-} "C";"L" <0pt>
\ar @{<-} "C";"R" <0pt>
 \endxy}
 \Ea
, 
   \Ba{c}\resizebox{44mm}{!}{  \xy
(-3,13)*+{_1}="UL";
(3,13)*+{_2}="UR";
(0,7)*{\bu}="U";
(0,-7)*+{_{2k}}="D";
(0,-10)*{\underbrace{\ \ \ \ \ \ \ \ \ \ \ \ \ \ \ \ \ }_{2k\ \text{valent vertex with no parallel edges}}};
    (0,0)*{\bu}="C";
  (7,0)*+{_3}="R";
   (-7,0)*+{_0}="L";
\ar @{<-} "U";"UL" <0pt>   
\ar @{<-} "U";"UR" <0pt>  
\ar @{<-} "C";"U" <0pt>
\ar @{<-} "C";"D" <0pt>
\ar @{->} "C";"L" <0pt>
\ar @{->} "C";"R" <0pt>
 \endxy}
 \Ea, \     \Ba{c}\resizebox{15mm}{!}{  \xy
 (0,14)*{_2}="1";
(0,-7)*{^0}="U";
(0,7)*{\bu}="D";
    (0,0)*{\bu}="C";
  (7,0)*{_3}="R";
   (-7,0)*{_1}="L";
\ar @{<-} "C";"U" <0pt>
\ar @{<-} "C";"D" <0pt>
\ar @{->} "C";"L" <0pt>
\ar @{->} "C";"R" <0pt>
\ar @{->} "D";"1" <0pt>
 \endxy}
 \Ea
 $$
Given any such a generator $\Ga\in Gr(\wh{\Delta}\PCYf_d^{marked})$, let $\widehat{\Ga}$ be the ribbon graph obtained from $\Ga$ by contracting all its operadic edges. As no operadic edges in $\Ga$ are connected to the $\caD\Lie_\infty$-type vertices of $\Ga$, the resulting hairy ribbon graph  $\widetilde{\Ga}$ can be identified with an element of the prop enveloping $P\caD\Lie_\infty$ of the dg properad $\caD\Lie_\infty$ to whose elements some number of ``static" corollas (the results of the above contractions) are attached;    the differential acts trivially on such ``static" corollas (hence their name); moreover we can assume without loss of generality  that vertices of these static corollas are distinguished, and all the half-edges attached to them are also distinguished. On the one hand, the complex
spanned by such ribbon graphs $\widehat{\Ga}$ computes, by Mashcke theorem, the cohomology of $Gr(\wh{\Delta}\PCYf_d^{marked})$, one the other hand it is isomorphic to the tensor products of trivial complexes (corresponding to ``static" corollas) with the prop closure of the properad  $\caD\Lie_\infty$
whose cohomology is generated \cite{L} by the ribbon corolla (\ref{1: DLie corolla}) modulo the relation (\ref{2: double Lie relations}). This result implies the isomorphism (\ref{2: final iso for wh(Delta)}) which in turn completes the proof of Claim II, and hence of the Main Theorem {\ref{1: Theorem on PCYf}}.
\end{proof}

\subsection{Proof of Theorem {\ref{1: Theorem on PCY}}}
In the previos subsection we proved that the epimorphism (see (\ref{2: map pi-marked}))
$$
\pi^{marked}: \PCYf_d^{marked}\lon \Delta \PCYf_d^{marked}
$$
is a quasi-isomorphism.
The above two complexes decompose into direct sums of subcomplexes
$$
\PCYf_d^{marked}=\bigoplus_{S\geq 0}   \PCYf_d^{marked}[S], \ \ \ \ \Delta \PCYf_d^{marked}=\bigoplus_{S\geq 0}  \Delta \PCYf_d^{marked}[S],
$$
spanned by hairy ribbon graphs with precisely $S$ marked sources, and the quasi-isomorphism $\pi^{marked}$ respects these decompositions. It remains to notice that 
$$
\PCYf_d^{marked}[0]=\PCY_d \ \ \text{and}\ \ \Delta \PCYf_d^{marked}[0]=
\Delta\PCY_d.
$$
The proof is completed.

\bip

\bip

{\Large
\section{\bf Deformation complex of properad $\PCYf_d$  and\\ moduli spaces $\cM_{g,1}$}
}

\sip

\subsection{Oriented version of the Kontsevich-Penner ribbon graph complex}
Let $\cM_{g,m}$ stand for the moduli space of genus $g$ algebraic curves with $m$ marked points. For any integer $d\in \Z$ there is a Kontsevich-Penner ribbon graph complex $\RGC_d$ which computes the compactly supported cohomology of these moduli spaces \cite{Ko0,Pe}
$$
H^\bu(\RGC_d)=\prod_{g\geq 0, m\geq 1\atop 2g+m\geq 3} H^{\bu +(d-1)m+d(2g-1)}_c(\cM_{g,m}).
$$
The complex $\RGC_d$ is spanned by at least trivalent ribbon graphs (with no hairs) which have with boundaries distinguished by integer labels, and whose edges  are directed up to a flip and multiplication by $(-1)^d$. The cohomological degree of such a generator $\Ga\in \RGC_d$ is defined by
$$
|\Ga|:= d(\# V(\Ga) -1) + (1-d)\#E(\Ga)
$$ 
where $V(\Ga)$ (resp. $E(\Ga)$) stands for the set of vertices (resp. edges)
of $\Ga$. The differential $\delta$ in $\RGC_d$ is given by the standard ``splitting of vertices" formula
\Beq\label{2: delta in RGC_0}
\delta\Ga= (-1)^{|\Ga|}\sum_{v\in V_\bu(\Ga)} \Ga\circ_v  \left(\xy
 (0,0)*{\bu}="a",
(5,0)*{\bu}="b",
\ar @{->} "a";"b" <0pt>
\endxy\right),
\Eeq
where 
$\Ga\circ_v  \left(\xy
 (0,0)*{\bu}="a",
(5,0)*{\bu}="b",
\ar @{->} "a";"b" <0pt>
\endxy\right)$ is a linear combination of ribbon graphs 
obtained from $\Ga$ by substituting into the vertex $v$ the ribbon graph $\xy
 (0,0)*{\bu}="a",
(5,0)*{\bu}="b",
\ar @{-} "a";"b" <0pt>
\endxy$ and then taking a sum over all possible reattachments of the edges (attached earlier to $v$) among the two newly created vertices in a way which respects the cyclic order of that edges and makes each new  vertex at least trivalent. In particular, if $v$ is itself trivalent, then $\Ga\circ_v  \left(\xy
 (0,0)*{\bu}="a",
(5,0)*{\bu}="b",
\ar @{-} "a";"b" <0pt>
\endxy\right)=0$.

\sip

There is a directed version $\dRGC_d$ of $\RGC_d$ in which the directions of edges are strictly fixed (i.e.\ no flipping is allowed), and the generating ribbon graphs are at least bivalent. There is a quasi-isomorphism of complexes
$$
\RGC_d \lon \dRGC_d
$$
which sends a ribbon graph $\Ga\in \RGC_d$  into a sum of graphs obtained from $\Ga$ by fixing directions on edges in all possible ways. Thus the complex $\dRGC_d$ gives us nothing new, but it contains a subcomplex $\ORGC_d\subset \dRGC_d$ spanned by directed ribbon graphs with no closed paths of directed edges; such ribbon graphs are often called {\it oriented}\, as there is a direction flow on such graphs starting at sources (that is, vertices with no incoming edges) and ending at targets (that is, vertices with no outgoing edges); such graphs are called sometimes ribbon {\it quivers}.

\sip

It was proven in \cite{Me2} that
$$
H^\bu(\ORGC_{d+1})\simeq  H^\bu(\RGC_d).
$$
Notice the change of the parameter $d$.
A second proof of this isomorphism of cohomology groups is given recently in \cite{Z}. The differentials in the complexes $\RGC_d$ and $\ORGC_{d+1}$
preserve  the genus $g$ and the numbers $m$ of boundaries  of its generators so that
they decompose into direct products of subcomplexes
$$
\RGC_{d}= \prod_{g\geq 0, m\geq 1\atop 2g+m\geq 3}\RGC_{d}^{g,m}, \ \ \ \ORGC_{d+1}= \prod_{g\geq 0, m\geq 1\atop 2g+m\geq 3}\ORGC_{d+1}^{g,m}.
$$
with 
\Beq\label{3: ORGCg,m and M_g,m}
H^\bu(\ORGC_{d+1}^{g,m})\simeq H^\bu(\RGC_{d}^{g,m})\simeq  H^{\bu +(d-1)m+d(2g-1)}_c(\cM_{g,m}).
\Eeq
The subspaces 
$$
\orgc_{d+1}:=\prod_{g\geq 1}\ORGC_{d+1}^{g,1}, \ \ \ \rgc_{d}:=\prod_{g\geq 1}\RGC_{d}^{g,1}
$$
spanned by ribbon graphs with precisely one boundary 
have the structures of dg Lie algebras \cite{MW}, and it was proven in \cite{Me2}
that  $\orgc_{d+1}$ and $\rgc_{d}$ are $\Lie_\infty$ quasi-isomorphic.
The Lie bracket in both Lie algebras $\orgc_{d+1}$ and $\rgc_{d}$ can be described by one and the same formula,
$$
[\Ga_1,\Ga_2]=\sum_{v\in V(\Ga_1)} \Ga_1\circ_v \Ga_2 -
(-1)^{|\Ga_1||\Ga_2|}\sum_{v\in V(\Ga_2)} \Ga_2\circ_v \Ga_1 ,
$$
where  $\Ga_1\circ_v \Ga_2$ stands for the linear combination of graphs obtain from $\Ga_1$ by substituting into its vertex $v$ the ribbon graph
$\Ga_2$
and taking the sum over all possible reattachments of the set edges attached earlier to $v$ to the set of corners of the unique boundary of $\Ga_2$ while respecting the cyclic structures of both sets (see \S 4.1 in \cite{MW} for more details).

\sip


\subsection{The derivation complex of $\PCYf_d$}   By a derivation complex  $\Der(\PCYf_d)$ of the dg properad $\PCYf_d$ we always understand (cf.\ \cite{MW2}) the dg Lie algebra of derivations of its genus completion $\wh{\PCYf_d}$.
As $\PCYf_d$ is a free properad, any derivation
$D \in \Der(\PCYf_d)$  is uniquely determined by its values on the $\bS$-bimodule  $E_d=\{E_d(m,n) \subset \PCYf(m,n)\}$ of the generators, so that we have an isomorphism of graded vector spaces
$$
 \Der(\PCYf_d)=
\prod_{m\geq 1,n\geq 0\atop m+n\geq 2} \Hom_{\bS_m^{op}\times \bS_n}\left(E_d(m,n), \widehat{\PCYf_d}(m,n)\right).
$$
As each $\bS_m^{op}\times \bS_n$-module $E_(m,n)$ is spanned by ribbon $(m,n)$-corollas of degree $(2-d)m -n+d$  (see (\ref{1: exs of ribb corollas}) for some examples of such corollas),
we can write
 \Beq\label{3: Der as product}
\Der(\PCYf_d)\simeq \prod_{m\geq 1,n\geq 0\atop m+n\geq 2} \widehat{\PCYf_d}(m\oplus n)[n-(2-d)m-d],
\Eeq
where $\widehat{\PCY_d}(m\oplus n)$ stands for the genus completed graded vector space of ribbon corollas having  $m+n$  hairs ($m$ of which are outgoing and $n$ are ingoing) which are {\it unlabelled, but instead are equipped with some fixed cyclic order}. For example, the differential $\delta$ in $\PCYf$
is a derivation  in $\Der(\PCYf_d)$
 given explicitly via the isomorphism (\ref{3: Der as product}) by the following infinite linear combination
\Beq\label{3: delta as a derivation}
\delta
=
\sum_{m\geq 1,n\geq 0\atop m+n\geq 2}    \underbrace{\resizebox{14mm}{!}{\xy
(-5,-6)*{}="a1",
(-5,6)*{}="a2",
(5,5)*{}="a3",
(5,-6)*{}="a4",
(8,0)*{}="a5",
(-8,0)*{}="a6",
(0,0)*+{\xy
 (0,0)*{\bu}="a",
(5,0)*{\bu}="b",
\ar @{->} "a";"b" <0pt>
\endxy}="b",
\ar @{<-} "a1";"b" <0pt>
\ar @{<-} "a2";"b" <0pt>
\ar @{<-} "a3";"b" <0pt>
\ar @{->} "a4";"b" <0pt>
\ar @{->} "a5";"b" <0pt>
\ar @{->} "a6";"b" <0pt>
\endxy}}_{\text{$m$ out-hairs, $n$ in-hairs}}
\Eeq
where each summand with fixed $m$ and $n$ stands for the sum  taken over all possible ways to attach the $m+n$ hairs 
to the two shown vertices while respecting their fixed  cyclic order.
 
 \sip

The linear space $\Der(\PCYf_d)$  has the standard Lie bracket denoted by $[\ ,\ ]$.
The differential $\delta$ in $\PCYf_d$ is an MC element of this Lie algebra, and hence it
makes  $\Der(\PCYf_d)$ into a {\it dg}\, Lie algebra equipped with the differential given by
$$
\p(D):=[\delta, D].
$$
We are interested in the associated cohomology $H^\bu(\Der(\PCYf_d),\p)$.
Given a derivation from $\Der(\PCYf_d)$ understood as a ribbon graph $\Ga\in \widehat{\PCYf_d}(m\oplus n)[n-(2-d)m-d]$ with cyclically ordered hairs, then
$$
\p\Ga=\delta\Ga + \sum_{m\geq 1,n\geq 0\atop m+n\geq 2}\hspace{-5mm}
\underbrace{\resizebox{12mm}{!}{
 \xy
(-5,-16)*{}="a1",
(-5,-4)*{}="a2",
(6,0)*+{\Ga}="a3",
(5,-16)*{}="a4",
(8,-10)*{}="a5",
(-8,-10)*{}="a6",
(0,-10)*{\bbu}="b",
\ar @{<-} "a1";"b" <0pt>
\ar @{<-} "a2";"b" <0pt>
\ar @{<-} "a3";"b" <0pt>
\ar @{->} "a4";"b" <0pt>
\ar @{->} "a5";"b" <0pt>
\ar @{->} "a6";"b" <0pt>
\endxy}}_{\text{ribbon $(m,n)$-corolla}}
\pm
\sum_{m,n\geq 1}\hspace{-5mm}
\overbrace{\resizebox{12mm}{!}{
 \xy
(-5,0)*{}="a1",
(-5,12)*{}="a2",
(6,12)*{}="a3",
(5,-4)*+{\Ga}="a4",
(8,6)*{}="a5",
(-8,6)*{}="a6",
(0,6)*{\bbu}="b",
\ar @{<-} "a1";"b" <0pt>
\ar @{<-} "a2";"b" <0pt>
\ar @{<-} "a3";"b" <0pt>
\ar @{->} "a4";"b" <0pt>
\ar @{->} "a5";"b" <0pt>
\ar @{->} "a6";"b" <0pt>
\endxy}}^{\text{ribbon $(m,n)$-corolla}}
$$
where $\delta$ is the standard ``splitting of vertices" differential in $\PCYf_d$ given by (\ref{2: delta in PreCYd}), and the second two sums are given by attaching $m,n)$-corollas to the in-hairs and, respectively, out-hairs corresponding hairs of $\Ga$. This formula implies that the differential $\p$ preserves the genus and the number of boundaries of the generators  $\Ga\in \widehat{\PCYf_d}(m\oplus n)[n-(2-d)m-d]$. Hence 
the derivation complex decomposes into the direct sum of subcomplexes
$$
\Der(\PCYf_d)=\bigoplus_{b\geq 1}\Der^{(b)}(\PCYf_d)
$$
parameterized by the number $b$ of boundaries of its generators. We are mostly interested in the subcomplex $\Der^{(1)}(\PCYf_d)$ spanned by derivations whose action on elements of $\PCYf_d$  preserve their number of boundaries. This subcomplex contains the differential $\delta$; it also contains the rescaling class
$$
r:=\sum_{m\geq 1,n\geq 0\atop m+n\geq 2}(m-n)\hspace{-5mm}  \underbrace{
\resizebox{10mm}{!}{\xy
(-5,-6)*{}="a1",
(-5,6)*{}="a2",
(5,5)*{}="a3",
(5,-6)*{}="a4",
(8,0)*{}="a5",
(-8,0)*{}="a6",
(0,0)*{\bbu}="b",
\ar @{<-} "a1";"b" <0pt>
\ar @{<-} "a2";"b" <0pt>
\ar @{<-} "a3";"b" <0pt>
\ar @{->} "a4";"b" <0pt>
\ar @{->} "a5";"b" <0pt>
\ar @{->} "a6";"b" <0pt>
\endxy}}_{\text{ribbon $(m,n)$-corolla}}
$$
corresponding to the derivation of $\PCYf_d$ which just rescales
$$
 C_{m,n} \lon \la^{m-n} C_{m,n}\ \ \ \forall \la\in \K^*,
$$
 the generating $(m,n)$-corollas. The derivation $r$ represents a non-trivial cohomology class 
in $H^0(\Der^{(1)}(\PCYf_d))$.

\sip

It is an important for us fact that the complex $(\Der(\PCYf_d), \p)$
can be identified (only as a complex, not as a dg Lie algbera) with the degree shifted deformation complex
\Beq\label{3: Der versus Def}
\Der(\PCYf_d) \simeq \Def(\PCYf_d \stackrel{\Id}{\rar} \PCYf_d)[1]=\bigoplus_{b\geq 1} \Def^{(b)}(\PCYf_d \stackrel{\Id}{\rar} \PCYf_d)[1]
\Eeq
of the identity automorphism of $\PCYf_d$ (we refer to \cite{MV} for the general deformation theory of morphisms of dg properads). The latter complex
also decomposes into the direct sum of subcomplexes.
$$
 \Def(\PCYf_d \rar \PCYf_d)[1]=\bigoplus_{b\geq 1} \Def^{(b)}(\PCYf_d \stackrel{\Id}{\rar} \PCYf_d)[1]
$$

\sip

The above discussion can be repeated for the dg properad $\PCY_d$ leading us to the complexes
$$
\Der(\PCY_d)=\bigoplus_{b\geq 1}\Der^{(b)}(\PCY_d)
\simeq \bigoplus_{b\geq 1} \Def^{(b)}(\PCY_d \rar \PCY_d)[1]
$$
which were studied earlier in \cite{Me2}.

\subsection{Proof of Theorem {\ref{1: Theorem on Def complex}}}
We want to show that $H^\bu(\Der^{(1)}(\PCYf_d))$ is highly non-trivial, that it contains all the cohomology groups 
 $$
H^\bu(\orgc_{d})=\prod_{g\geq 1} H_c^{\bu-1+2g(d-1)}(\cM_{g,1})\ \ \oplus \ \ \K\langle r \rangle
$$ 
of the moduli spaces $\cM_{g,1}$ (which in turn span, up to a rescaling class $r$,  the cohomology of the above mentioned complex $\orgc_{d}$ of oriented ribbon graphs with precisely one boundary. Here $\K\langle r \rangle$ is the 1-dimensional vector space spanned by the rescaling class $r$.

\sip

There is an almost obvious monomorphism of complexes (cf.\ \cite{MW2})
\Beq\label{3: Morhism F from orgc to Der}
\Ba{rccl}
 f: & \orgc_d &\to & \Der^{(1)}(\PCYf_d)\simeq\Def^{(1)}(\PCYf_d \rar \PCYf_d)[1]\\
         &   \ga & \to & f(\ga)
         \Ea
\Eeq
where the derivation $f(\ga)$ acts (from the right) on the generating ribbon
$(m,n)$-corolla of $\PCYf_d$ by attaching  $m+n$ hairs to the corners of the
 unique boundary $b$ of $\ga$ in all possible ways while respecting the cyclic 
 orders of both sets,
 $$
 f(\Ga): \underbrace{
\resizebox{10mm}{!}{\xy
(-5,-6)*{}="a1",
(-5,6)*{}="a2",
(5,5)*{}="a3",
(5,-6)*{}="a4",
(8,0)*{}="a5",
(-8,0)*{}="a6",
(0,0)*{\bbu}="b",
\ar @{<-} "a1";"b" <0pt>
\ar @{<-} "a2";"b" <0pt>
\ar @{<-} "a3";"b" <0pt>
\ar @{->} "a4";"b" <0pt>
\ar @{->} "a5";"b" <0pt>
\ar @{->} "a6";"b" <0pt>
\endxy}}_{\text{ribbon $(p,q)$-corolla}}
\ \ \ \lon \ 
 \sum_{\text{maps}\ [p+q]\rar C(b)} \resizebox{14mm}{!}{\xy
(-5,-6)*{}="a1",
(-5,6)*{}="a2",
(5,5)*{}="a3",
(5,-6)*{}="a4",
(8,0)*{}="a5",
(-8,0)*{}="a6",
(0,0)*+{\ga}="b",
\ar @{<-} "a1";"b" <0pt>
\ar @{<-} "a2";"b" <0pt>
\ar @{<-} "a3";"b" <0pt>
\ar @{->} "a4";"b" <0pt>
\ar @{->} "a5";"b" <0pt>
\ar @{->} "a6";"b" <0pt>
\endxy},
 $$
A similar morphism of complexes 
$$
 f^3:  \orgc_d \to  \Def^{(1)}(\PCY_d \rar \PCY_d)[1]\
$$
holds true for the second dg properad. The canonical projection $ep: \PCYf_d\rar \PCY_d$ (which is just the quotient by the ideal generated by $(m,0)$-corollas) induces 
a morphism of complexes
$$
e:  \Def^{(1)}(\PCYf_d \stackrel{\Id}{\rar} \PCYf_d)[1]\lon  \Def^{(1)}(\PCYf_d \stackrel{e}{\rar} \PCYf_d)[1].
$$
As the map $e$ must vanish on the generating corollas of $\PCYf_d$ with no incoming hairs, we conclude that there is an identification 
$$ 
\Def^{(1)}(\PCYf_d \stackrel{e}{\rar} \PCY_d)[1]= \Def^{(1)}(\PCY_d \stackrel{\Id}{\rar} \PCY_d)[1]
$$
of complexes. Hence we obtain a commutative diagram
$$
\xymatrix{
 \orgc_d \ar[dr]_-{f^3}\ar[r]^{f\ \ \ \ \ \ \ \ \ \ \ \ \ \ \ } &  \Def^{(1)}(\PCYf_d \stackrel{\Id}{\rar} \PCYf_d)[1] \ar[d]^-{e}\\
 & \Def^{(1)}(\PCY_d \stackrel{\Id}{\rar} \PCY_d)[1]
}
$$
It was proven in \S 5.8 of  \cite{Me2} that the morphism $f^3$ is a quasi-isomorphism. Hence the above diagram implies a monomorphism
of the cohomology groups
\Beq\label{3: mono H(orgc) to Def(PCYf)}
0\lon H^\bu(\orgc_d) \lon H^{\bu+1}(\Def^{(1)}(\PCYf_d \stackrel{\Id}{\rar} \PCYf_d))
\Eeq
which in turn implies the required Theorem~{\ref{1: Theorem on Def complex}}. $\Box$

\subsection{Proof of Corollary {\ref{1: Corollary on TZ}}}
Assume the contrary, i.e.\ assume that the properad $\cP\cV^{(d)}$ is Koszul. In this case the canonical epimorphism from the cobar construction on the Koszul dual cooperad to $\cP\cV^{(d)}$,
$$
s: \PCYf_d=\Omega((\cP\cV^{(1-d)})^{\ac}) \lon \cP\cV^{(1-d)},
$$
is a quasi-isomophism. The map $s$ induces a morphism of the deformation complexes (by composing the identity map with the projection $s$)
$$
s: \Def^{(1)}(\PCYf_d \stackrel{\Id}{\rar} \PCYf_d)\lon 
\Def^{(1)}(\PCYf_d \stackrel{s }{\rar} \cP\cV^{(1-d)})
$$
which must also be a quasi-isomorphism by the general theory \cite{MV}.
Then the monomorphism (\ref{3: mono H(orgc) to Def(PCYf)}) implies the injection of cohomology groups
\Beq\label{3: injection of H(orgc) into H(derPCYf)}
0\lon H^\bu(\orgc_d) \lon H^{\bu+1}(\Def^{(1)}(\PCYf_d \stackrel{s}{\rar}  \cP\cV^{(1-d)})).
\Eeq
Since the properad $\cP\cV^{(1-d)}$ is generated by only 2- and 3-valent corollas, the ribbon graphs having at least one 4-valent vertex belong to the kernel of the composition
$$
s\circ f: \orgc_d\lon \Def^{(1)}(\PCYf_d \stackrel{\Id}{\rar} \PCYf_d)\lon \Def^{(1)}(\PCYf_d \stackrel{s}{\rar}  \cP\cV^{(1-d)}).
$$
Hence the injection (\ref{3: injection of H(orgc) into H(derPCYf)})  implies that every cohomology class in the totality $\prod_{g\geq 1} H_c^{\bu-1+2g(d-1)}(\cM_{g,1})$ can be represented by ribbon graphs in $\orgc_d$ having at most  trivalent vertices which is impossible (see e.g. \S 2.11 in \cite{Me2} or \S 5.10.1 in \cite{Me1})). This contradiction proves Corollary {\ref{1: Corollary on TZ}}. $\Box$

\def\cprime{$'$}


\begin{thebibliography}{10}

\bibitem[A]{A} Marcelo Aguiar, {\em Infinitesimal Hopf Algebras}, In: New trends in Hopf algebra theory,
vol.\ 267, Contemporary Math. (2000), 1-29.


\bibitem[GK]{GK}
Victor~Ginzburg and Michail Kapranov.
\newblock {\it Koszul duality for operads},
\newblock {Duke Math. J.}, 76(1):203--272, 1994.



\bibitem[IKV]{IKV} Natalia Iyudu, Maxim Kontsevich, and Yannis Vlassopoulos. {\em Pre-CalabiYau
algebras as noncommutative Poisson structures},  Journal of Algebra {\bf 567} (2021), pp. 63-90.


\bibitem[Kh]{Kh} Anton Khoroshkin, {\it Input/output coloring and Gröbner basis for dioperads}, arXiv:2602.20069 (2026).

\bibitem[Ko]{Ko0} Maxim Kontsevich, {\em Intersection Theory on the Moduli Space of Curves
and the Matrix Airy Function}, Commun.\ Math.\ Phys. {\bf 147} (1992)  1-23


\bibitem[KTV]{KTV}
Maxim Kontsevich, Alex Takeda, and Yiannis Vlassopoulos {\em  Pre-Calabi-Yau algebras
and topological quantum field theories},  Eur. J. Math.,  {\bf 11} (2025), article number 15.


\bibitem[L]{L} Johan Leray, {\it  Protoperads II: Koszul duality}, Journal de l'\'{E}cole polytechnique — Math\'{e}matiques, {\bf 7} (2020) 897-941.

\bibitem[LV]{LV} Johan Leray and Bruno Vallette, {\em Pre-Calabi-Yau algebras and homotopy
double Poisson gebras}, arXiv: 2203.05062 (2022).




\bibitem[M1]{Me1} Sergei Merkulov,  {\it A complex of ribbon quivers and $\cM_{g,m}$}, arXiv:2503.02020 (2025)

\bibitem[M2]{Me2} Sergei Merkulov,  {\it A low-valence ribbon graph complex computing the cohomology of $\cM_{g,m}$},  arXiv:2605.04950 (2026)

\bibitem[MV]{MV}  Sergei Merkulov and  Bruno Vallette,
{\em Deformation theory of representations of prop(erad)s I \& II},
 J.\ f\"ur die reine und angewandte Mathematik (Qrelle)  {\bf 634}, 51-106,
 \& {\bf 636}, 123-174 (2009)



 \bibitem[MW1]{MW} Sergei Merkulov and Thomas Willwacher, {\em Props of 
 involutive Lie bialgebras and moduli spaces of curves},  arXiv:1511.07808
  (2015) 51pp.
  
  \bibitem[MW2]{MW2} Sergei Merkulov and Thomas Willwacher, {\it Deformation theory of
 Lie bialgebra properads},   In: Geometry and Physics: A Festschrift in honour of
 Nigel Hitchin, Oxford University Press 2018, pp. 219-248





\bibitem[P]{Pe}  Robert Penner, {\em The decorated Teichm\"uller space of punctured surfaces}, Comm.\ Math.\
Phys. {\bf 113} (1987), 299-339.

\bibitem[PT]{PT} Kate Poirier and Thomas Tradler,
{\it Koszuality of the $V^{(d)}$ dioperad}, J. Homotopy Relat. Struct., 
{\bf  14} (2019)  477-507.




\bibitem[Q]{Q} Alexandre Quesney, {\it Balanced infinitesimal bialgebras, 
double Poisson gebras and pre-Calabi-Yau algebras}, arXiv:2312.14893
(2023)

\bibitem[TZ]{TZ} Thomas Tradler and Mahmoud Zeinalian, {\em  Algebraic string operations},
K-Theory, {\bf 38}(1) (2007) 59-82.


\bibitem[VdB]{vdB} 
Michel Van den Bergh, {\it Double Poisson algebras}, Trans.\ Amer.\ Math.\ Soc. {\bf 360} (2008) 5711-5769.



\bibitem[Z]{Z} Marko\ \v Zivkovi\' c, {\em Direct quasi-isomorphism from ribbon graph complex to oriented ribbon graph complex},  preprint  2026.     

 \end{thebibliography}
\end{document}